\documentclass[12pt,a4paper]{amsart}

\usepackage{euscript,amsfonts,amssymb,amsmath,amscd,amsthm,enumerate,hyperref,mathrsfs}

\usepackage{tikz,bm,float}
\usetikzlibrary{calc}
\colorlet{lgray}{white!85!black}
\colorlet{lred}{white!75!red}
\newcommand{\bra}[1]{\left\langle #1\right|}
\newcommand{\ket}[1]{\left|#1\right\rangle}

\usetikzlibrary{decorations.pathreplacing}

\usepackage{comment}

\usepackage{graphicx}

\usepackage{color}

\usepackage[margin=1.1in]{geometry}
\newtheorem{theorem}{Theorem} 
\newtheorem*{theorem*}{Theorem}
\newtheorem{lemma}[theorem]{Lemma}
\newtheorem{definition}[theorem]{Definition}
\newtheorem{proposition}[theorem]{Proposition}

\newtheorem{corollary}[theorem]{Corollary}

\theoremstyle{remark}
\newtheorem{remark}[theorem]{Remark}
\newtheorem{example}[theorem]{Example}

\numberwithin{equation}{section} \numberwithin{theorem}{section}

\newcommand{\braket}[2]{\langle #1 | #2 \rangle}

\newcommand{\Pf}{\ensuremath{\mathrm{Pf}}}
\newcommand{\kernel}{\mathsf{K}}
\newcommand{\Ikernel}{\mathsf{I}}
\newcommand{\Rkernel}{\mathsf{R}}
\newcommand{\sign}{\ensuremath{\mathrm{sign}}}
\newcommand{\I}{\ensuremath{\mathbf{i}}}
\def \Ai {{\rm Ai}}

\usepackage{MnSymbol}

\usepackage{skak}

\usepackage{adjustbox}
\usetikzlibrary{arrows}
\usetikzlibrary{decorations.pathreplacing}

\begin{document}

\title{The second-class particle in the half-line open TASEP}

\author[Kailun Chen]{Kailun Chen}\address{Institute of Mathematics, Leipzig University, Augustusplatz 10, 04109 Leipzig, Germany.}\email{kailunxsx@gmail.com}

\begin{abstract}
We consider the second-class particle in half-line open TASEP under two different initial conditions with shock discontinuities. The exact formulas for the distribution of the second-class particle can be derived by using the color--position symmetry theorem of colored half-space TASEP. We study the asymptotic distributions of second-class particles under the constant scaling and KPZ-type scaling.
\end{abstract}

\maketitle

\section{Introduction}

\subsection{Overview}
The totally asymmetric simple exclusion process (TASEP) is a fundamental stochastic transport model, with applications in several areas of out-of-equilibrium physics, including biological micro-systems and traffic flow. 
Notably, it establishes intriguing theoretical connections with Kardar-Parisi-Zhang universality \cite{Cor12} and random matrix theory \cite{Fer10}. 
Because of its simple definition and analytical tractability, TASEP has played an important role in the study of non-equilibrium statistical mechanics. 
The model is defined on a one-dimensional lattice, where each site can contain at most one particle. 
Each particle attempts to jump to a neighboring empty site, and the term ``totally asymmetric'' refers to the fact that jumps are allowed only in one direction. 
Over time, many variants of the model have been introduced, incorporating additional features motivated by different physical settings.

A first important issue is the introduction of the totally asymmetric exclusion process on the positive integers with a single particle source at the origin (we call it half-line open TASEP). Half-line open TASEP has been of interest due to the rich phase diagram that is expected to hold for limiting distributions based on the boundary strength. The work by \cite{Ligg75} reveals a pivotal insight: the long-term behavior of this process undergoes a phase transition. If the particle production rate at the source and the initial density fall below a critical threshold, the stationary measure is a product measure; otherwise, it exhibits spatial correlation. Building upon the methodology presented in \cite{DEHP93}, \cite{Gro04} demonstrated that these correlations can be elucidated through a matrix product representation. This realization has yielded various outcomes within the physics literature, encompassing the derivation of hydrodynamic limits, comprehension of phase diagrams, and the establishment of large deviation principles. A substantial body of literature on this subject exists, exemplified by \cite{DMZ18} and related references. The asymptotic analysis of the half-line open TASEP with a step initial condition has been explored in \cite{BBCS18a} (in the equivalent LPP language, also referencing prior works \cite{BR01a, BR01b, SI04}) and \cite{BBCS18b} (in the equivalent FTASEP language). These studies unveil that starting from an empty configuration, the total number of particles in the system at time $t$ exhibits Gaussian fluctuations on the scale of $t^{1/2}$ when the particle production rate at the source is strictly below 1/2. However, when the production rate surpasses 1/2, Tracy-Widom GSE fluctuations \cite{TW96} manifest on the $t^{1/3}$ scale. In the critical case, where the production rate equals 1/2, Tracy-Widom GOE fluctuations \cite{TW96} also emerge on the $t^{1/3}$ scale. A two-dimensional crossover among the GUE \cite{TW94}, GOE, and GSE distributions was elucidated through multipoint distribution analysis. The asymptotic analysis of the half-line open TASEP with a stationary initial condition has been conducted in \cite{BFO20} and \cite{BFO22} (in an equivalent LPP language). These works reveal analogues of Baik-Rains distribution \cite{BR00} and the Airy stat process of Baik-Ferrari-Péché \cite{BFP10} in the half-space geometry. Notably, recent advancements have significantly advanced our understanding of the half-line open ASEP, a one-parameter generalization of the half-line open TASEP. For further insights, one may refer to \cite{BBCW18, BBC20, IMS22, He23}.

Among other interesting issues, the second-class particles have been frequently considered in the literature. The second-class particles play a crucial role, particularly in the presence of shock discontinuities, where their position can effectively indicate the location of the shock (see \cite{Ligg99}, Chapter 3). Consequently, several investigations have been dedicated to understanding second-class particles. In the scenario of Bernoulli-Bernoulli initial conditions, it has been demonstrated that the fluctuations of second-class particles exhibit Gaussian behavior on a scale of $t^{1/2}$, and the limiting process corresponds to a Brownian motion (see \cite{Ferr92, FF94}). However, when dealing with shocks having non-random initial conditions, the situation differs due to the emergence of typical Kardar-Parisi-Zhang (KPZ) fluctuations from the dynamics. For initial conditions featuring non-random densities $\lambda < \rho$, \cite{FGN19} has determined that the fluctuations of the second-class particle, starting from the origin, follow the difference of two independent random variables with GOE Tracy-Widom distributions on a $t^{1/3}$ scale. In \cite{FN23}, the limit process is characterized by the difference of two independent $\mathrm{Airy}_1$ processes. More classes of initial conditions were explored in \cite{BF22}, where identities for the one-point distribution of a second-class particle were derived by employing a symmetry theorem within a multi-color TASEP framework \cite{BB21}. Additionally, asymptotic fluctuations were analyzed under various scaling conditions. For more recent works on shocks in interacting particle systems, see e.g. \cite{FN15,FN17,Nej18,FN20,Nej20} and the references therein.

In this paper, we consider the second-class particles in the half-line open TASEP. Following \cite{BF22}, we will consider two initial conditions (with one shock or two shocks) for half-line open TASEP which can be explored by employing the color--position symmetry property of colored half-space TASEP. The main results of this paper are that we write down the exact distribution of the second-class particles in terms of the height function of standard half-line open TASEP with step initial condition, and then the asymptotic results of second-class particles can be obtained by using the properties of the latter. We introduce our models and explain the asymptotic results below.

\subsection{Multi-color half-line open TASEP with one shock or two shocks}
\label{ssec:models}
In this paper, we will introduce several multi-color TASEPs on $\mathbb{N}$ with open boundary. All of them are continuous-time processes with second-class particles which start from distinct initial configurations.

\begin{figure}[h]
\begin{tikzpicture}[scale=1.25]
\node[circle, inner sep=10pt, draw=black, fill=orange, opacity=0.7] (a) at (-3.3, 0) {};
\draw (-3.3, 0) node[scale=0.6,blue] {$\text{reservoir}$};
\node[circle, inner sep=3.8pt, draw=black, fill=white, opacity=0.7] (b) at (-2, 0) {};
\node[circle, inner sep=3.8pt, draw=black, fill=white, opacity=0.7] (c) at (-1, 0) {};
\node[circle, inner sep=3.8pt, draw=black, fill=white, opacity=0.7] (d) at (0, 0) {};
\node[circle, inner sep=3.8pt, draw=black, fill=red, opacity=0.7] (e) at (1, 0) {};
\node at (1, 0) {$2$};
\node[circle, inner sep=3.8pt, draw=black, fill=orange, opacity=0.7] (f) at (2, 0) {};
\node at (2, 0) {$1$};
\node[circle, inner sep=3.8pt, draw=black, fill=orange, opacity=0.7] (g) at (3, 0) {};
\node at (3, 0) {$1$};
\node[circle, inner sep=3.8pt, draw=black, fill=orange, opacity=0.7] (h) at (4, 0) {};
\node at (4, 0) {$1$};
\node[circle, inner sep=3.8pt, draw=black, fill=white, opacity=0.7] (i) at (5, 0) {};
\node[circle, inner sep=3.8pt, draw=black, fill=white, opacity=0.7] (j) at (6, 0) {};
\node[circle, inner sep=3.8pt, draw=black, fill=white, opacity=0.7] (k) at (7, 0) {};
%
\draw (a) -- (b);
\draw (b) -- (c);
\draw (c) -- (d);
\draw (d) -- (e);
\draw (e) -- (f);
\draw (f) -- (g);
\draw (g) -- (h);
\draw (h) -- (i);
\draw (i) -- (j);
\draw (j) -- (k);
\draw[dotted, thick] (k) -- (8,0);
\draw (-2, -0.5) node[scale=0.5] {$1$};
\draw (-1, -0.5) node[scale=0.5] {$2$};
\draw (-0.5, -0.5) node[scale=0.5] {$\cdots$};
\draw (0, -0.5) node[scale=0.5] {$M_1$};
\draw (1, -0.5) node[scale=0.5] {$M_1+1$};
\draw (2, -0.5) node[scale=0.5] {$M_1+2$};
\draw (3, -0.5) node[scale=0.5] {$M_1+3$};
\draw (3.5, -0.5) node[scale=0.5] {$\cdots$};
\draw (4.1, -0.5) node[scale=0.5] {$M_1+M_2+1$};
\draw[-stealth, dashed, semithick] (-3.3,0.4) to[out=45,in=135] node[midway, above] {$\alpha$} (-2,0.2) ;
\draw[-stealth, dashed, semithick] (1,0.2) to[out=45,in=135] node[midway] {$\times$} (2,0.2) ;
\draw[-stealth, dashed, semithick] (4,0.2) to[out=45,in=135] node[midway, above] {$1$} (5,0.2) ;
\begin{scope}[shift={(0,-1.8)}]
\node[circle, inner sep=10pt, draw=black, fill=orange, opacity=0.7] (a) at (-3.3, 0) {};
\draw (-3.3, 0) node[scale=0.6,blue] {$\text{reservoir}$};
\node[circle, inner sep=3.8pt, draw=black, fill=white, opacity=0.7] (b) at (-2, 0) {};
\node[circle, inner sep=3.8pt, draw=black, fill=white, opacity=0.7] (c) at (-1, 0) {};
\node[circle, inner sep=3.8pt, draw=black, fill=white, opacity=0.7] (d) at (0, 0) {};
\node[circle, inner sep=3.8pt, draw=black, fill=red, opacity=0.7] (e) at (1, 0) {};
\node at (1, 0) {$2$};
\node[circle, inner sep=3.8pt, draw=black, fill=red, opacity=0.7] (f) at (2, 0) {};
\node at (2, 0) {$2$};
\node[circle, inner sep=3.8pt, draw=black, fill=red, opacity=0.7] (g) at (3, 0) {};
\node at (3, 0) {$2$};
\node[circle, inner sep=3.8pt, draw=black, fill=red, opacity=0.7] (h) at (4, 0) {};
\node at (4, 0) {$2$};
\node[circle, inner sep=3.8pt, draw=black, fill=white, opacity=0.7] (i) at (5, 0) {};
\node[circle, inner sep=3.8pt, draw=black, fill=white, opacity=0.7] (j) at (6, 0) {};
\node[circle, inner sep=3.8pt, draw=black, fill=white, opacity=0.7] (k) at (7, 0) {};
%
\draw (a) -- (b);
\draw (b) -- (c);
\draw (c) -- (d);
\draw (d) -- (e);
\draw (e) -- (f);
\draw (f) -- (g);
\draw (g) -- (h);
\draw (h) -- (i);
\draw (i) -- (j);
\draw (j) -- (k);
\draw[dotted, thick] (k) -- (8,0);
\draw (-2, -0.5) node[scale=0.5] {$1$};
\draw (-1, -0.5) node[scale=0.5] {$2$};
\draw (-0.5, -0.5) node[scale=0.5] {$\cdots$};
\draw (0, -0.5) node[scale=0.5] {$M_1$};
\draw (1, -0.5) node[scale=0.5] {$M_1+1$};
\draw (2, -0.5) node[scale=0.5] {$M_1+2$};
\draw (3, -0.5) node[scale=0.5] {$M_1+3$};
\draw (3.5, -0.5) node[scale=0.5] {$\cdots$};
\draw (4.1, -0.5) node[scale=0.5] {$M_1+M_2+1$};
\draw[-stealth, dashed, semithick] (-3.3,0.4) to[out=45,in=135] node[midway, above] {$\alpha$} (-2,0.2) ;
\draw[-stealth, dashed, semithick] (4,0.2) to[out=45,in=135] node[midway, above] {$1$} (5,0.2) ;
\end{scope}
\end{tikzpicture}
\caption{The half-line open TASEP with one-shock initial condition. Upper: half-line open TASEP $\eta_{t}$ with initial configuration $\eta_{0}$; Bottom: half-line open TASEP $\tilde{\eta}_{t}$ with initial configuration $\tilde{\eta}_{0}$. }\label{fig:one-shock}
\end{figure}
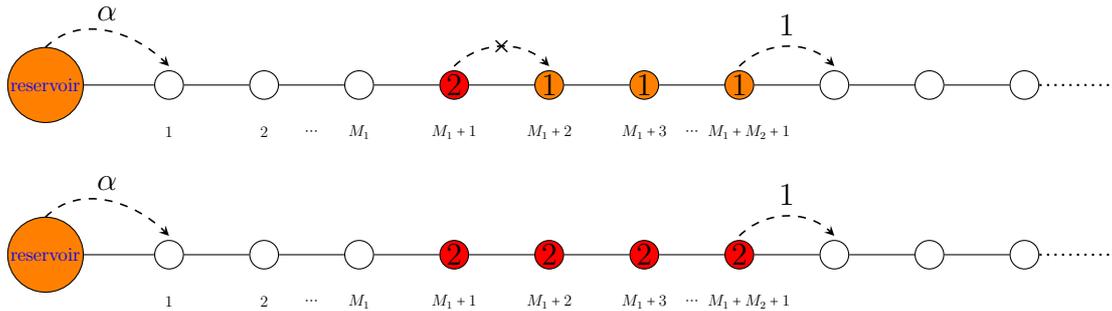

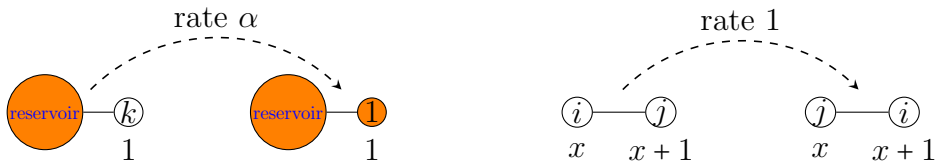
\begin{figure}[h]
\begin{tikzpicture}
\node[circle, inner sep=10pt, draw=black, fill=orange, opacity=0.7] (aL) at (-1.1, 0) {};
\draw (-1.1, 0) node[scale=0.6,blue] {$\text{reservoir}$};

\node[circle, inner sep=3.8pt, draw=black] (bL) at (0, 0) {};
\node at (0, 0) {$k$};
\node at (0, -0.5) {$1$};

\draw (aL) -- (bL);

\node[circle, inner sep=10pt, draw=black, fill=orange, opacity=0.7] (bR) at (1.1+1, 0) {};
\draw (1.1+1, 0) node[scale=0.6,blue] {$\text{reservoir}$};

\node[circle, inner sep=3.8pt, draw=black, fill=orange, opacity=0.7] (aR) at (2.2+1, 0) {};
\node at (2.2+1, 0) {$1$};
\node at (2.2+1, -0.5) {$1$};

\draw (aR) -- (bR);

\draw[-stealth, dashed, semithick] (-0.5,0.3) to[out=45,in=135] node[midway, above] {rate $\alpha$} (2.8,0.3) ;
\begin{scope}[shift={(7,0)}]
\node[circle, inner sep=4pt, draw=black] (aL) at (-1.1, 0) {};
\node at (-1.1, 0) {$i$};
\node at (-1.1, -0.5) {$x$};

\node[circle, inner sep=4pt, draw=black] (bL) at (0, 0) {};
\node at (0, 0) {$j$};
\node at (0, -0.5) {$x+1$};

\draw (aL) -- (bL);

\node[circle, inner sep=4pt, draw=black] (bR) at (1.1+1, 0) {};
\node at (1.1+1, 0) {$j$};
\node at (1.1+1, -0.5) {$x$};

\node[circle, inner sep=4pt, draw=black] (aR) at (2.2+1, 0) {};
\node at (2.2+1, 0) {$i$};
\node at (2.2+1, -0.5) {$x+1$};

\draw (aR) -- (bR);

\draw[-stealth, dashed, semithick] (-0.5,0.3) to[out=45,in=135] node[midway, above] {rate $1$} (2.6,0.3) ;
\end{scope}
\end{tikzpicture}
\caption{The transition rates, where $k \in \{2, +\infty\}$, $i < j$ and $i,j \in \{1,2, +\infty\}$.}\label{fig:rates}
\end{figure}

The first half-line open TASEP has three colors (we denote the first-class particle by 1, second-class particle by 2 and hole by $+\infty$), which is defined on the state-space parametrized by occupation variables $\{(\eta(x))_{x\in\mathbb{Z}_{> 0}} \in \{1,2,+\infty\}^{\mathbb{Z}_{> 0}}\}$. The initial configuration $\eta_{0}(x)$ has one second-class particle which is located between $M_1$ holes (on the  interval $[1, M_1]$) and $M_2$ first-class particles (on the interval $[M_1+2, M_1+M_2+1]$), the sites larger than $M_1+M_2+1$ are all filled with holes. See the upper panel of Figure \ref{fig:one-shock}. The state $\eta_{t}(x)$ at given time $t\in[0,\infty)$ and position $x\in\mathbb{Z}_{> 0}$ evolves according to the following dynamics: a particle $\eta_{t}(x)$ jumps from site $x$ to $x+1$ at exponential rate 1 if $\eta_{t}(x)<\eta_{t}(x+1)$, a first-class particle is created at site 1 at rate $\alpha$ if $\eta_{t}(1)>1$. Otherwise, we keep the state unchanged. See Figure \ref{fig:rates}. Let $f(t)$ represent the position of the unique second-class particle in the process $\eta_{t}$.

The second half-line open TASEP (we denote by $\tilde{\eta}_{t}$) is similar to the first one, but with a different initial condition. The second initial configuration (which we denote by $\tilde{\eta}_{0}$) can be obtained from the first one by replacing the first-class particles on the interval $[M_1+2, M_1+M_2+1]$ with the second-class particles. See the lower panel of Figure \ref{fig:one-shock}. Obviously, the initial conditions $\eta_0$ and $\tilde{\eta}_0$ create one shock. Let $\mathbf{N}_{c}(x,t)$ denote the number of particles of color $c$ in the process $\tilde{\eta}_t$ at time $t$ that are located weakly to the right of $x\in \mathbb{R}_{\geq 0}$, where $c=1,2$.

\begin{figure}[h]
\begin{tikzpicture}[scale=1.3]
\node[circle, inner sep=10pt, draw=black, fill=orange, opacity=0.7] (a) at (-3.3, 0) {};
\draw (-3.3, 0) node[scale=0.6,blue] {$\text{reservoir}$};
\node[circle, inner sep=3.8pt, draw=black, fill=white, opacity=0.7] (b) at (-2, 0) {};
\node[circle, inner sep=3.8pt, draw=black, fill=white, opacity=0.7] (c) at (-1, 0) {};
\node[circle, inner sep=3.8pt, draw=black, fill=orange, opacity=0.7] (d) at (0, 0) {};
\node at (0, 0) {$1$};
\node[circle, inner sep=3.8pt, draw=black, fill=orange, opacity=0.7] (e) at (1, 0) {};
\node at (1, 0) {$1$};
\node[circle, inner sep=3.8pt, draw=black, fill=white, opacity=0.7] (f) at (2, 0) {};
\node[circle, inner sep=3.8pt, draw=black, fill=white, opacity=0.7] (g) at (3, 0) {};
\node[circle, inner sep=3.8pt, draw=black, fill=red, opacity=0.7] (h) at (4, 0) {};
\node at (4, 0) {$2$};
\node[circle, inner sep=3.8pt, draw=black, fill=orange, opacity=0.7] (i) at (5, 0) {};
\node at (5, 0) {$1$};
\node[circle, inner sep=3.8pt, draw=black, fill=orange, opacity=0.7] (j) at (6, 0) {};
\node at (6, 0) {$1$};
\node[circle, inner sep=3.8pt, draw=black, fill=white, opacity=0.7] (k) at (7, 0) {};
%
\draw (a) -- (b);
\draw (b) -- (c);
\draw (c) -- (d);
\draw (d) -- (e);
\draw (e) -- (f);
\draw (f) -- (g);
\draw (g) -- (h);
\draw (h) -- (i);
\draw (i) -- (j);
\draw (j) -- (k);
\draw[dotted, thick] (k) -- (8,0);
\draw (-2, -0.5) node[scale=0.48] {$1$};
\draw (-1.5, -0.5) node[scale=0.48] {$\cdots$};
\draw (-1, -0.5) node[scale=0.48] {$N$};
\draw (0, -0.5) node[scale=0.48] {$N+1$};
\draw (0.5, -0.5) node[scale=0.48] {$\cdots$};
\draw (1, -0.5) node[scale=0.48] {$N+M$};
\draw (2, -0.5) node[scale=0.48] {$N+M+1$};
\draw (2.5, -0.5) node[scale=0.48] {$\cdots$};
\draw (3, -0.5) node[scale=0.48] {$N+2M$};
\draw (4, -0.5) node[scale=0.48] {$N+2M+1$};
\draw (4.9, -0.5) node[scale=0.48] {$N+2M+2$};
\draw (5.5, -0.5) node[scale=0.48] {$\cdots$};
\draw (6.1, -0.5) node[scale=0.48] {$2N+2M+1$};
\draw[-stealth, dashed, semithick] (-3.3,0.4) to[out=45,in=135] node[midway, above] {$\alpha$} (-2,0.2) ;
\draw[-stealth, dashed, semithick] (1,0.2) to[out=45,in=135] node[midway, above] {$1$} (2,0.2) ;
\draw[-stealth, dashed, semithick] (4,0.2) to[out=45,in=135] node[midway] {$\times$} (5,0.2) ;
\begin{scope}[shift={(0,-1.8)}]
\node[circle, inner sep=10pt, draw=black, fill=orange, opacity=0.7] (a) at (-3.3, 0) {};
\draw (-3.3, 0) node[scale=0.6,blue] {$\text{reservoir}$};
\node[circle, inner sep=3.8pt, draw=black, fill=white, opacity=0.7] (b) at (-2, 0) {};
\node[circle, inner sep=3.8pt, draw=black, fill=white, opacity=0.7] (c) at (-1, 0) {};
\node[circle, inner sep=3.8pt, draw=black, fill=red, opacity=0.7] (d) at (0, 0) {};
\node at (0, 0) {$2$};
\node[circle, inner sep=3.8pt, draw=black, fill=red, opacity=0.7] (e) at (1, 0) {};
\node at (1, 0) {$2$};
\node[circle, inner sep=3.8pt, draw=black, fill=white, opacity=0.7] (f) at (2, 0) {};
\node[circle, inner sep=3.8pt, draw=black, fill=white, opacity=0.7] (g) at (3, 0) {};
\node[circle, inner sep=3.8pt, draw=black, fill=green, opacity=0.7] (h) at (4, 0) {};
\node at (4, 0) {$3$};
\node[circle, inner sep=3.8pt, draw=black, fill=green, opacity=0.7] (i) at (5, 0) {};
\node at (5, 0) {$3$};
\node[circle, inner sep=3.8pt, draw=black, fill=green, opacity=0.7] (j) at (6, 0) {};
\node at (6, 0) {$3$};
\node[circle, inner sep=3.8pt, draw=black, fill=white, opacity=0.7] (k) at (7, 0) {};
%
\draw (a) -- (b);
\draw (b) -- (c);
\draw (c) -- (d);
\draw (d) -- (e);
\draw (e) -- (f);
\draw (f) -- (g);
\draw (g) -- (h);
\draw (h) -- (i);
\draw (i) -- (j);
\draw (j) -- (k);
\draw[dotted, thick] (k) -- (8,0);
\draw (-2, -0.5) node[scale=0.48] {$1$};
\draw (-1.5, -0.5) node[scale=0.48] {$\cdots$};
\draw (-1, -0.5) node[scale=0.48] {$N$};
\draw (0, -0.5) node[scale=0.48] {$N+1$};
\draw (0.5, -0.5) node[scale=0.48] {$\cdots$};
\draw (1, -0.5) node[scale=0.48] {$N+M$};
\draw (2, -0.5) node[scale=0.48] {$N+M+1$};
\draw (2.5, -0.5) node[scale=0.48] {$\cdots$};
\draw (3, -0.5) node[scale=0.48] {$N+2M$};
\draw (4, -0.5) node[scale=0.48] {$N+2M+1$};
\draw (4.9, -0.5) node[scale=0.48] {$N+2M+2$};
\draw (5.5, -0.5) node[scale=0.48] {$\cdots$};
\draw (6.1, -0.5) node[scale=0.48] {$2N+2M+1$};
\draw[-stealth, dashed, semithick] (-3.3,0.4) to[out=45,in=135] node[midway, above] {$\alpha$} (-2,0.2) ;
\draw[-stealth, dashed, semithick] (1,0.2) to[out=45,in=135] node[midway, above] {$1$} (2,0.2) ;
\draw[-stealth, dashed, semithick] (6,0.2) to[out=45,in=135] node[midway, above] {$1$} (7,0.2) ;
\end{scope}
\end{tikzpicture}
\caption{The half-line open TASEP with two-shock initial condition. Upper: half-line open TASEP $\xi_{t}$ with initial configuration $\xi_{0}$; Bottom: half-line open TASEP $\tilde{\xi}_{t}$ with initial configuration $\tilde{\xi}_{0}$. }\label{fig:two-shocks}
\end{figure}
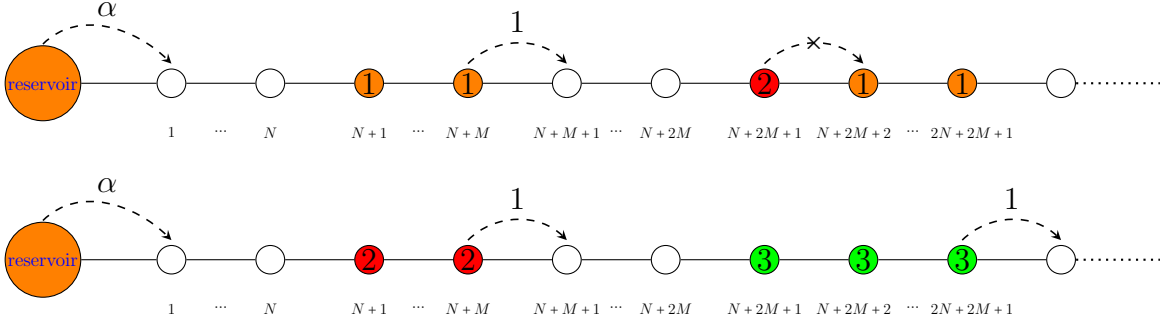

We extend the one-shock initial configuration to the two-shock case in the following way. In the two-shock case of the first initial configuration $\eta_0$ (we denote by $\xi_0$), we put $N$ holes on the interval $[1, N]$, $M$ first-class particles on the interval $[N+1, N+M]$, $M$ holes on the interval $[N+M+1, N+2M]$, one second-class particle at the position $N+2M+1$, $N$ first-class particles on the interval $[N+2M+2, 2N+2M+1]$, and the sites larger than $2N+2M+1$ are all filled with holes. We denote the half-line open TASEP with initial condition $\xi_0$ by $\xi_t$. See the upper panel of Figure \ref{fig:two-shocks}. The dynamics of $\xi_t$ are also given by Figure \ref{fig:rates}. We denote the position of the unique second-class particle in the process $\xi_{t}$ by $g(t)$. 

The two-shock case of the second initial configuration $\tilde{\eta}_0$ (we denote by $\tilde{\xi}_0$) can be obtained from $\xi_0$ by replacing the particles on the interval $[N+1, N+M]$ with the second-class particles, and replacing the particles on the interval $[N+2M+1, 2N+2M+1]$ with the third-class particles (which we denote by 3). We denote the half-line open TASEP with initial condition $\tilde{\xi}_0$ by $\tilde{\xi}_t$. See the lower panel of Figure \ref{fig:two-shocks}. The dynamics of $\tilde{\xi}_t$ are also given by Figure \ref{fig:rates}, but with $k \in \{2, 3, +\infty\}$, $i < j$ and $i,j \in \{1,2, 3, +\infty\}$. Let $\mathsf{N}_{i}(x,t)$ be the number of particles of color $i$ in $\tilde\xi_{t}$ which are weakly to the right of $x\in \mathbb{R}$, where $i \in \{1,2,3\}$.

\subsection{Asymptotic distribution of the second-class particles under the constant scaling}
\label{ssec:constant-scaling}
In this section, we consider the asymptotic probability distribution of the second-class particles when $M_1,M_2$ are fixed constants. Before stating our main results, let us introduce some notation. For any $L\in\mathbb{N}$ and any state $\eta=(\eta_1,\cdots,\eta_L) \in \{0,1\}^{L}$, we denote $\ell(\eta)=\#\{i: \eta_i=1\}$ and $A(\eta)=\#\{\text{clusters of 1's in}\ \eta \}$. For $1 \leq j \leq A(\eta)$, we denote the number of 1 in the $j$-th cluster of 1's by $\sigma_j(\eta)$. It is easy to see that $\sum_{j=1}^{A(\eta)}\sigma_{j}(\eta)=\ell(\eta)$. For $1 \leq j \leq A(\eta)-1$, we denote the number of 0 in the cluster of 0's between the $j$-th and $(j+1)$-th cluster of 1's by $\tau_j(\eta)$. And $\tau_0(\eta)$ is the number of 0 before the first cluster of 1's, $\tau_{A(\eta)}(\eta)$ is the number of 0 after the last cluster of 1's. Then we have $\sum_{k=0}^{A(\eta)}\tau_{k}(\eta)=L-\ell(\eta)$. We will use the notations $\Psi_i(\eta)=\sum_{j=i}^{A(\eta)}\sigma_{j}(\eta)$ for $1 \leq i \leq A(\eta)+1$, and $\Phi_j(\eta)= \sum_{k=j}^{A(\eta)}\tau_{k}(\eta)$ for $0 \leq j \leq A(\eta)+1$, where $\Psi_{A(\eta)+1}(\eta)=\Phi_{A(\eta)+1}(\eta)=0$. For example, if we take $L=12$ and $\eta=(0,0,1,0,0,1,1,0,1,1,0,0)$, then we have $\ell(\eta)=5$, $A(\eta)=3$ and 
\[\sigma_1(\eta)=1, \sigma_2(\eta)=2, \sigma_3(\eta)=2;\]
\[\tau_0(\eta)=2,\tau_1(\eta)=2, \tau_2(\eta)=1,\tau_3(\eta)=2;\]
\[\Psi_1(\eta)=5, \Psi_2(\eta)=4, \Psi_3(\eta)=2, \Psi_4(\eta)=0;\]
\[\Phi_0(\eta)=7, \Phi_1(\eta)=5,\Phi_2(\eta)=3,\Phi_3(\eta)=2,\Phi_4(\eta)=0.\]

\begin{figure}
\begin{center}
\begin{tikzpicture}[scale=1]
\draw [blue] (2,12) node {$(5,7)$};
\draw[black,line width=1pt,->] (1.9,11.8) -- (1,11.3) node[midway, anchor=west]{$c/\alpha$};
\draw [blue] (1,11) node {$(5,6)$};
\draw[black,line width=1pt,->] (0.9,10.8) -- (0,10.3) node[midway, anchor=west]{$c/\alpha$};
\draw [black] (0,10) node {$(5,5)$};
\draw[black,line width=1pt,->] (-0.1,9.8) -- (-1,9.3) node[midway, anchor=west]{$c$};
\draw (-1,9) node {$(5,4)$};
\draw[black,line width=1pt,->] (0.1,9.8) -- (1,9.3) node[midway, anchor=east]{$c$};
\draw [blue] (1,9) node {$(4,5)$};
\draw[black,line width=1pt,->] (-1.1,8.8) -- (-2,8.3) node[midway, anchor=west]{$c$};
\draw [black] (-2,8) node {$(5,3)$};
\draw[black,line width=1pt,->] (-0.9,8.8) -- (0,8.3) node[midway, anchor=east]{$c$};
\draw [blue] (0,8) node {$(4,4)$};
\draw[black,line width=1pt,->] (-2.1,7.8) -- (-3,7.3) node[midway, anchor=west]{$c$};
\draw [black] (-3,7) node {$(5,2)$};
\draw[black,line width=1pt,->] (-1.9,7.8) -- (-1,7.3) node[midway, anchor=east]{$c$};
\draw [black] (-1,7) node {$(4,3)$};
\draw[black,line width=1pt,->] (-3.1,6.8) -- (-4,6.3) node[midway, anchor=west]{$c$};
\draw (-4,6) node {$(5,1)$};
\draw[black,line width=1pt,->] (-2.9,6.8) -- (-2,6.3) node[midway, anchor=east]{$c$};
\draw [black] (-2,6) node {$(4,2)$};
\draw[black,line width=1pt,->] (-1.1,6.8) -- (-2,6.3) node[midway, anchor=west]{$c$};
\draw[black,line width=1pt,->] (-0.9,6.8) -- (0,6.3) node[midway, anchor=east]{$c$};
\draw (0,6) node {$(3,3)$};
\draw[black,line width=1pt,->] (-4.1,5.8) -- (-5,5.3) node[midway, anchor=west]{$c$};
\draw[black,line width=1pt,->] (-3.9,5.8) -- (-3,5.3) node[midway, anchor=east]{$c$};
\draw[black,line width=1pt,->] (-2.1,5.8) -- (-3,5.3) node[midway, anchor=west]{$c$};
\draw[black,line width=1pt,->] (-1.9,5.8) -- (-1,5.3) node[midway, anchor=east]{$c$};
\draw[black,line width=1pt,->] (-0.1,5.8) -- (-1,5.3) node[midway, anchor=west]{$c$};
\draw[black,line width=1pt,->] (0.1,5.8) -- (0.9,5.3) node[midway, anchor=east]{$c$};
\draw [black] (-5,5) node {$(5,0)$};
\draw (-3,5) node {$(4,1)$};
\draw [black] (-1,5) node {$(3,2)$};
\draw [blue] (1,5) node {$(2,3)$};
\draw[black,line width=1pt,->] (-3.1,4.8) -- (-4,4.3) node[midway, anchor=west]{$c$};
\draw[black,line width=1pt,->] (-2.9,4.8) -- (-2,4.3) node[midway, anchor=east]{$c$};
\draw[black,line width=1pt,->] (-1.1,4.8) -- (-2,4.3) node[midway, anchor=west]{$c$};
\draw[black,line width=1pt,->] (-0.9,4.8) -- (0,4.3) node[midway, anchor=east]{$c$};
\draw [black] (-4,4) node {$(4,0)$};
\draw (-2,4) node {$(3,1)$};
\draw [black] (0,4) node {$(2,2)$};
\draw[black,line width=1pt,->] (-2.1,3.8) -- (-3,3.3) node[midway, anchor=west]{$c$};
\draw[black,line width=1pt,->] (-1.9,3.8) -- (-1,3.3) node[midway, anchor=east]{$c$};
\draw[black,line width=1pt,->] (-0.1,3.8) -- (-1,3.3) node[midway, anchor=west]{$c$};
\draw[black,line width=1pt,->] (0.1,3.8) -- (1,3.3) node[midway, anchor=east]{$c$};
\draw [black] (-3,3) node {$(3,0)$};
\draw (-1,3) node {$(2,1)$};
\draw (1,3) node {$(1,2)$};
\draw[black,line width=1pt,->] (-1.1,2.8) -- (-2,2.3) node[midway, anchor=west]{$c$};
\draw[black,line width=1pt,->] (-0.9,2.8) -- (0,2.3) node[midway, anchor=east]{$c$};
\draw[black,line width=1pt,->] (0.9,2.8) -- (0,2.3) node[midway, anchor=west]{$c$};
\draw[black,line width=1pt,->] (1.1,2.8) -- (2,2.3) node[midway, anchor=east]{$c$};
\draw [black] (-2,2) node {$(2,0)$};
\draw (0,2) node {$(1,1)$};
\draw [blue] (2,2) node {$(0,2)$};
\draw[black,line width=1pt,->] (-0.1,1.8) -- (-1,1.3) node[midway, anchor=west]{$c$};
\draw[black,line width=1pt,->] (0.1,1.8) -- (1,1.3) node[midway, anchor=east]{$c$};
\draw [black] (-1,1) node {$(1,0)$};
\draw [blue] (1,1) node {$(0,1)$};
\draw[black,line width=1pt,->] (0.9,8.8) -- (0.1,8.3) node[midway, anchor=west]{$c/\alpha$};
\draw[black,line width=1pt,->] (-0.1,7.8) -- (-0.9,7.3) node[midway, anchor=west]{$c/\alpha$};
\draw[black,line width=1pt,->] (0.9,4.8) -- (0.1,4.3) node[midway, anchor=west]{$c/\alpha$};
\draw[black,line width=1pt,->] (1.9,1.8) -- (1,1.3) node[midway, anchor=west]{$c/\alpha$};
\draw[black,line width=1pt,->] (0.9,0.8) -- (0,0.3) node[midway, anchor=west]{$c/\alpha$};
\draw [black] (0,0) node {$(0,0)$};
\end{tikzpicture}
\end{center}
\caption{The DEHP-tree of $\eta=(0,0,1,0,0,1,1,0,1,1,0,0)$; the DEHP-tree of $\eta$ starts from the node $(5,7)$ and ends with the nodes $(k,0)$ for $0 \leq k \leq 5$; the nodes $(5,7)$, $(5,6)$, $(4,5)$, $(4,4)$, $(2,3)$, $(0,2)$, $(0,1)$ have only one child, the other nodes have two children.}
\label{fig:DHEP-tree}
\end{figure}
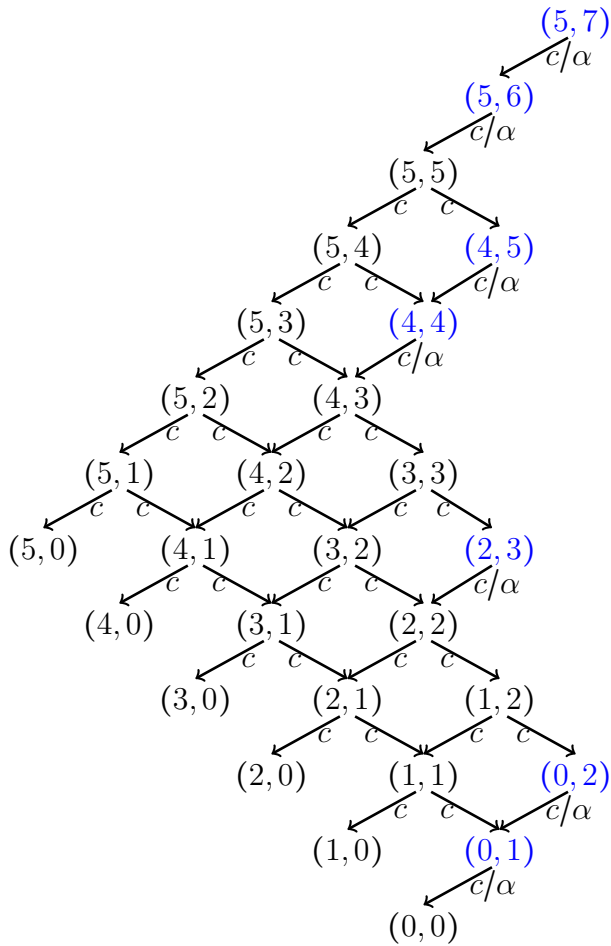

\begin{definition}
\label{def:DEHP-tree}
Let $\alpha \in (0,1)$ and $c=\frac14$. For any $L\in\mathbb{N}$ and any state $\eta=(\eta_1,\cdots,\eta_L) \in \{0,1\}^{L}$, we define the DEHP-tree of $\eta$ in the following way:
\begin{itemize}
\item start from a node represented by a tuple $(\Psi_1(\eta), \Phi_0(\eta))$, which we call the {\em root}.
\item any node can have either one child or two children. If a node $(x,y)$ has two children, we generate two new nodes: the left one $(x,y-1)$ and the right one $(x-1,y)$, and draw a {\em weighted directed path} $\swarrow$ (resp. $\searrow$) pointing from the {\em parent} $(x,y)$ to the {\em left child} (resp. the {\em right child}), where each path has {\em weight} $c$. If a node $(x,y)$ has only one child, we generate one new node: the left child $(x,y-1)$, and the parent $(x,y)$ is connected with the left child by a path $\swarrow$ weighted by $c/\alpha$. 
\item When $(x,y)\in\{(x,y): x=\Psi_{i}(\eta)\ \text{and}\ \Phi_i(\eta) < y \leq \Phi_{i-1}(\eta) \ \text{for}\ 1 \leq i \leq  A(\eta)+1\}$, the node $(x,y)$ has only one child. When $(x,y)\in\{(x,y): x=\Psi_{i}(\eta)\ \text{and}\ 1 \leq y \leq \Phi_{i}(\eta) \ \text{for}\ 1 \leq i \leq  A(\eta)\} \cup \{(x,y): \Psi_{i+1}(\eta) <x<\Psi_{i}(\eta), 1 \leq y \leq \Phi_{i}(\eta)  \ \text{for}\ 1 \leq i \leq  A(\eta)\}$, the node $(x,y)$ has two children. 
\item end with the nodes $(k,0)$ where $k=0,1,\cdots,\Psi_{1}(\eta)$ which we call the {\em endpoints}.
\end{itemize}
For every node $(x,y)$ reached by this construction, there is at least one directed path from the root $(\Psi_1(\eta),\Phi_0(\eta))$ to $(x,y)$.
We define the weight of such a path as the product of the weights along the path. We define the partition function of a node $(x,y)$ as the sum of weights of all possible directed paths from the root $(\Psi_1(\eta), \Phi_0(\eta))$ to that node, and denote it by $Z_{x,y}(\eta)$. In particular, we denote the partition function of endpoints $(k,0)$ as $Z_{k}(\eta)$, where $k=0,1,\cdots,\Psi_{1}(\eta)$. See Figure \ref{fig:DHEP-tree} for illustration. 
\end{definition}

Our first main result is the asymptotic behavior of the second-class particles in $\eta_t$ and $\tilde{\eta}_{t}$ when $\alpha\leq \frac12$:
\begin{theorem}[Theorem \ref{thm:one-shock-stationary-4.1}]
Let $\alpha\leq \frac12$. The asymptotic distribution  of the unique second-class particle existing in $\eta_t$ is given by:
\begin{align*}
\lim_{t\to\infty}\mathbb{P}(f(t)>0)=\sum_{n=0}^{M_1}\binom{M_1+M_2+1}{n}\alpha^n(1-\alpha)^{M_1+M_2+1-n}.
\end{align*}

Let $0 \leq m \leq M_2$, the asymptotic distribution  of the height function of second-class particles in $\tilde{\eta}_{t}$ is given by:
\begin{align*}
\lim_{t\to\infty}\mathbb{P}(\mathbf{N}_{2} (1,t)=m)=\binom{M_1+M_2+1}{M_1+M_2+1-m}\alpha^{M_1+M_2+1-m}(1-\alpha)^{m}.
\end{align*}
\end{theorem}

Our second main result is the asymptotic behavior of the second-class particles in $\eta_t$ and $\tilde{\eta}_{t}$ when $\alpha > \frac12$:

\begin{theorem}[Theorem \ref{thm:one-shock-stationary-4.2}]
Let $\alpha > \frac12$ and we use the notation $B_{n}=\{\eta \in \{0,1\}^{M_1+M_2+1}: \ell(\eta)=n\}$. The asymptotic distribution of the unique second-class particle existing in $\eta_t$ is given by:
\begin{align*}
\lim_{t \to \infty } \mathbb{P}(f(t)>0)=
\sum_{n=0}^{M_1}\sum_{\eta \in B_{n}}\sum_{k=0}^{\Psi_{1}(\eta)} \frac{Z_k(\eta)}{2^k} \cdot \left(1+\frac{\alpha-\frac12}{\alpha}k\right).
\end{align*}

Let $0 \leq m \leq M_2$, the asymptotic distribution  of the height function of second-class particles in $\tilde{\eta}_{t}$ is given by:
\begin{align*}
\lim_{t\to\infty}\mathbb{P}(\mathbf{N}_{2} (1,t)=m)
=\sum_{\eta \in B_{M_1+M_2+1-m}}\sum_{k=0}^{\Psi_{1}(\eta)} \frac{Z_k(\eta)}{2^k} \cdot \left(1+\frac{\alpha-\frac12}{\alpha}k\right).
\end{align*}
\end{theorem}

Our third main result is the asymptotic behavior of the second-class particle in $\xi_t$ and the third-class particles in $\tilde{\xi}_{t}$ when $ \alpha \leq \frac12$:
\begin{theorem}[Theorem \ref{thm:two-shocks-stationary-4.1}]
Let $ \alpha \leq \frac12$. The asymptotic distribution of the unique second-class particle existing in the half-line open TASEP $\xi_t$ is given by:
\begin{multline*}
\lim_{t\to\infty}\mathbb{P}(g(t)>0)=\sum_{m=0}^{N-1}\sum_{n=0}^{M+N-m}\binom{M+N}{m}\binom{M+N+1}{n}\alpha^{m+n}(1-\alpha)^{2M+2N+1-m-n}\\
+\sum_{m=N}^{M+N}\sum_{n=0}^{M}\binom{M+N}{m}\binom{M+N+1}{n}\alpha^{m+n}(1-\alpha)^{2M+2N+1-m-n}.
\end{multline*}

Let $0 \leq s \leq N$, the asymptotic distribution of the height function of third-class particles in $\tilde{\xi}_t$ is given by:
\begin{multline*}
\lim_{t\to\infty}\mathbb{P}(\mathsf{N}_{3} (1,t)=s)=\sum_{n=1}^{s}\binom{M+N}{N-n}\binom{M+N+1}{M+N+1-s+n}\alpha^{M+2N+1-s}(1-\alpha)^{M+s}\\
+\sum_{m=N}^{M+N}\binom{M+N}{m}\binom{M+N+1}{M+N+1-s}\alpha^{M+N+1+m-s}(1-\alpha)^{M+N-m+s}.
\end{multline*}
\end{theorem}

Our fourth main result is the asymptotic behavior of the second-class particle in $\xi_t$ and the third-class particles in $\tilde{\xi}_{t}$ when $ \alpha > \frac12$:
\begin{theorem}[Theorem \ref{thm:two-shocks-stationary-4.2}]
Let $ \alpha > \frac12$ and we use the notation $B_{m,n}=\{\eta \in \{0,1\}^{2M+2N+1}: \sum_{i=1}^{M+N}\eta_i=m, \sum_{j=M+N+1}^{2M+2N+1}\eta_j=n\}$. The asymptotic distribution of the unique second-class particle existing in the half-line open TASEP $\xi_t$ is given by:
\begin{multline*}
\lim_{t\to\infty}\mathbb{P}(g(t)>0)=\sum_{m=0}^{N-1}\sum_{n=0}^{M+N-m}\sum_{\eta \in B_{m,n}} \sum_{k=0}^{\Psi_{1}(\eta)} \frac{Z_k(\eta)}{2^k} \cdot \left(1+\frac{\alpha-\frac12}{\alpha}k\right)  \\
+\sum_{m=N}^{M+N}\sum_{n=0}^{M}\sum_{\eta \in B_{m,n}}\sum_{k=0}^{\Psi_{1}(\eta)} \frac{Z_k(\eta)}{2^k} \cdot \left(1+\frac{\alpha-\frac12}{\alpha}k\right).
\end{multline*}

Let $0 \leq s \leq N$, the asymptotic distribution of the height function of third-class particles in $\tilde{\xi}_t$ is given by:
\begin{multline*}
\lim_{t\to\infty}\mathbb{P}(\mathsf{N}_{3} (1,t)=s)=\sum_{n=1}^{s}\sum_{\eta \in B_{N-n,M+N+1-s+n}} \sum_{k=0}^{\Psi_{1}(\eta)} \frac{Z_k(\eta)}{2^k} \cdot \left(1+\frac{\alpha-\frac12}{\alpha}k\right)  \\
+\sum_{m=N}^{M+N}\sum_{\eta \in B_{m,M+N+1-s}}\sum_{k=0}^{\Psi_{1}(\eta)} \frac{Z_k(\eta)}{2^k} \cdot \left(1+\frac{\alpha-\frac12}{\alpha}k\right).
\end{multline*}
\end{theorem}

\begin{remark}
The asymptotic distribution of the height function of second-class particles in $\tilde{\xi}_t$ is similar to \eqref{eq:one-shock-height-1} and \eqref{eq:one-shock-height-2}. For the
sake of brevity, we only provide the asymptotic distribution of the height function of third-class particles in $\tilde{\xi}_t$.
\end{remark}

\subsection{Asymptotic distribution of the second-class particles under the KPZ scaling}
In this section, we consider the asymptotic probability distribution of the second-class particles under the KPZ scaling $t^{\frac23}$. Below we will use the \emph{half-space ${Airy}_2$ process} $\mathcal{A}_{\varpi}^{HS}$ which was introduced in \cite{BBCS18a}, see Section \ref{sec:KPZ-scaling} for the definition. Our first main result is the asymptotic behavior of the second-class particles in $\eta_t$ and $\tilde{\eta}_{t}$:
\begin{theorem}[Theorem \ref{thm:one-shock-KPZ-5}]
\label{thm:one-shock-KPZ}
Assume that $M_1=M_2=\lfloor at^{\frac23} \rfloor$, $a\in \mathbb{R}_{>0}$. Set $\alpha = \frac{1+2^{4/3}\varpi t^{-1/3}}{2}$, where $\varpi \in \mathbb{R}$. Then the asymptotic distribution  of the unique second-class particle disappearing from the half-line open TASEP $\eta_t$ is given by:
\begin{align}
\label{eq:one-shock-disappear}
\lim_{t \to \infty}\mathbb{P}(f(t)<0)  = \mathbb{P}\left(\mathcal{A}_{\varpi}^{HS}(0)+2^{4/3}a^2 \leq \mathcal{A}_{\varpi}^{HS}(2^{2/3}a)\right).
\end{align}
The asymptotic distribution  of the height function of second-class particles in $\tilde{\eta}_{t}$ is given by:
\begin{align}
\label{one-shock-kpz}
\lim_{t\to\infty}\frac{\mathbf{N}_{2} (1,t)-at^{2/3}}{t^{1/3}} \stackrel{d}{=} \min \left\{2^{-4/3}\mathcal{A}_{\varpi}^{HS}(0)+a^2 - 2^{-4/3}\mathcal{A}_{\varpi}^{HS}(2^{2/3}a),0 \right\}.
\end{align}
\end{theorem}

Our second main result is the asymptotic behavior of the second-class particle in $\xi_t$ and the third-class particles in $\tilde{\xi}_{t}$:
\begin{theorem}[Theorem \ref{thm:two-shocks-KPZ-5}]
\label{thm:two-shocks-KPZ-5}
Assume that $M=N=\lfloor at^{\frac23} \rfloor$, $a\in \mathbb{R}_{>0}$. Set $\alpha = \frac{1+2^{4/3}\varpi t^{-1/3}}{2}$, where $\varpi \in \mathbb{R}$. Then the asymptotic distribution of the unique second-class particle disappearing from the half-line open TASEP $\xi_t$ is given by:
\begin{multline}
\label{eq:two-shocks-disappear}
\lim_{t \to \infty}\mathbb{P}(g(t)<0)  = \mathbb{P}\left(2^{-\frac{4}{3}}\mathcal{A}_{\varpi}^{HS}(2^{\frac53}a)-2^{-\frac{4}{3}}\mathcal{A}_{\varpi}^{HS}(2^{\frac23}a) - 3a^2\right.\\
\left. + \min\left\{2^{-\frac{4}{3}}\mathcal{A}_{\varpi}^{HS}(2^{\frac23}a)-2^{-\frac{4}{3}}\mathcal{A}_{\varpi}^{HS}(0) - a^2,0 \right\} \geq 0\right).
\end{multline}
The asymptotic distribution of the height function of third-class particles in $\tilde{\xi}_t$ is given by:
\begin{multline}
\label{eq:two-shocks-kpz}
\lim_{t\to\infty}\frac{\mathsf{N}_{3} (1,t)-at^{2/3}}{t^{1/3}} \stackrel{d}{=} \min \left\{0, 2^{-4/3}\mathcal{A}_{\varpi}^{HS}(2^{2/3}a)+3a^2 - 2^{-4/3}\mathcal{A}_{\varpi}^{HS}(2^{5/3}a) \right.\\
\left. +\max\{0, 2^{-4/3}\mathcal{A}_{\varpi}^{HS}(0)+a^2 - 2^{-4/3}\mathcal{A}_{\varpi}^{HS}(2^{2/3}a)\}\right\}.
\end{multline}
\end{theorem}

\subsection{Outline}
The paper is organized as follows. In Section \ref{sec:ColorPositionSymmetry} we introduce the colored half-space TASEP and the color--position symmetry theorem. By using the color--position symmetry theorem, we derive the exact distribution of the second-class particles in different kinds of multi-color half-line TASEPs in Section \ref{sec:distributional-identity}. In Section \ref{sec:constant-scaling} and Section \ref{sec:KPZ-scaling}, we write down the asymptotic results for the second-class particles under the constant-scaling and KPZ-scaling.

\subsection*{Acknowledgments} 
I am grateful to A. Bufetov for very useful discussions. The work was partially supported by A. Bufetov's Grant from the European Research Council (ERC), Grant Agreement No. 101041499.

\section{Color--position symmetry theorem of colored half-space TASEP}
\label{sec:ColorPositionSymmetry}

In this section, we first define the colored half-space TASEP. Subsequently, we describe the evolution of colored half-space TASEP as the random walk on the type $B/C$ Hecke algebra, and introduce the algebraic properties of the colored half-space TASEP. Finally, we will use the algebraic structure of colored half-space TASEP to derive the color--position symmetry theorem. 

\subsection{Colored half-space TASEP}
\label{ssec:Chs-TASEP}

The \emph{colored half-space totally asymmetric simple exclusion process} (colored half-space TASEP) is described as a continuous-time Markov process operating on the set of bijections from $\mathbb{N} \cup -\mathbb{N}$ to itself. Initially, particles of various colors are placed at corresponding positions, with each particle at position $x$ assigned the color $x$. Formally, the configurations of the system are represented by functions $\eta$, where $\eta(x)$ denotes the color of the particle at position $x$. The process commences with the identity function, $\eta(x) = x$.

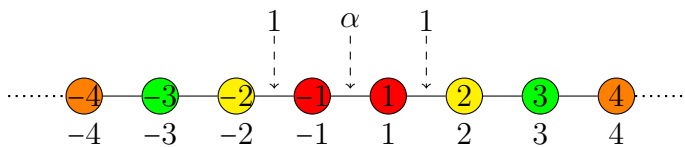
\begin{figure}[h]
\begin{tikzpicture}
\node[circle, inner sep=4.8pt, draw=black, fill=red, opacity=0.7] (a) at (0, 0) {};
\node at (0, 0) {$-1$};
\node[circle, inner sep=4.8pt, draw=black, fill=red, opacity=0.7] (b) at (1, 0) {};
\node at (1, 0) {$1$};
\node[circle, inner sep=4.8pt, draw=black, fill=yellow, opacity=0.7] (c) at (-1, 0) {};
\node at (-1, 0) {$-2$};
\node[circle, inner sep=4.8pt, draw=black, fill=yellow, opacity=0.7] (d) at (2, 0) {};
\node at (2, 0) {$2$};
\node[circle, inner sep=4.8pt, draw=black, fill=green, opacity=0.7] (e) at (-2, 0) {};
\node at (-2, 0) {$-3$};
\node[circle, inner sep=4.8pt, draw=black, fill=green, opacity=0.7] (f) at (3, 0) {};
\node at (3, 0) {$3$};
\node[circle, inner sep=4.8pt, draw=black, fill=orange, opacity=0.7] (g) at (-3, 0) {};
\node at (-3, 0) {$-4$};
\node[circle, inner sep=4.8pt, draw=black, fill=orange, opacity=0.7] (h) at (4, 0) {};
\node at (4, 0) {$4$};
\draw (a) -- (b);
\draw (b) -- (d);
\draw (d) -- (f);
\draw (f) -- (h);
\draw (c) -- (a);
\draw (e) -- (c);
\draw (g) -- (e);
\draw[dotted, thick] (-4,0) -- (g);
\draw[dotted, thick] (h) -- (5,0);
\node at (0, -0.5) {$-1$};
\node at (1, -0.5) {$1$};
\node at (2, -0.5) {$2$};
\node at (3, -0.5) {$3$};
\node at (4, -0.5) {$4$};
\node at (-1, -0.5) {$-2$};
\node at (-2, -0.5) {$-3$};
\node at (-3, -0.5) {$-4$};
\draw[black, dashed, ->] (0.5,0.8)--(0.5,0.1) ;
\node at (0.5, 1) {$\alpha$};
\draw[black, dashed, ->] (1.5,0.8)--(1.5,0.1) ;
\node at (1.5, 1) {$1$};
\draw[black, dashed, ->] (-0.5,0.8)--(-0.5,0.1) ;
\node at (-0.5, 1) {$1$};
\end{tikzpicture}
\caption{The colored half-space TASEP with packed initial condition $\eta(x)=x$; the transition rate associated with edge $(-1,1)$ is $\alpha$; for $x \geq 1$, the transition rates associated with edge $(x, x+1)$ and $(-x-1, -x)$ are given by $1$.}\label{fig:colored-half space-TASEP}
\end{figure}

The dynamics of this model involve associating exponential clocks of rate 1 with each edge $(x, x+1)$ for $x \geq 1$, and an exponential clock of rate $\alpha$ with the edge $(-1, 1)$. When a clock attached to the edge $(x,x+1)$ with rate 1 rings, particles are swapped at positions $x$ and $x+1$ if $\eta(x)<\eta(x+1)$, and a similar swap occurs for particles at positions $-x$ and $-x-1$ at the same time. When the clock with rate $\alpha$ rings, a swap occurs between $\eta(-1)$ and $\eta(1)$ if $\eta(-1)<\eta(1)$. Otherwise, we leave the system configuration unchanged. See Figure \ref{fig:colored-half space-TASEP} for an illustration. 
For any fixed time $t$, almost surely one can find a sequence $i_k\in\mathbb N$ such that no clock associated with the edge $(i_k,i_k+1)$ rings up to time $t$.
Given this sequence, the process can be sampled within the finite intervals $[i_k+1, i_{k+1}]$, ensuring the well-defined nature of this process (see the graphical construction in \cite{Har78}).

\subsection{Random walk on the type $B/C$ Hecke algebra}
\label{ssec:RW-Hecke}

Consider the group $B_n$, which is the hyperoctahedral group of signed permutations (referred to as permutations) on $n$ letters. 
In other words, $B_n$ consists of permutations $\pi$ of ${-n,\cdots, -1, 1, \cdots, n}$ such that $\pi(-i)=-\pi(i)$. 
We treat this group as a Coxeter group, having generators $s_0=(-1,1)$ and $s_k=(k,k+1)$ for $0<k<n$. 
The \emph{length} $l(\pi)$ is defined as the minimum number of generators required to express $\pi$. 
We can express $l(\pi)$ as the sum $l_{0}(\pi)+l_{1}(\pi)$, where $l_{0}$ counts the occurrences of $s_0$ and $l_{1}$ counts the occurrences of $s_k$ for $k>0$ in any decomposition $\pi=s_{k_1} \cdots s_{k_n}$ with $n$ being minimal. 
This decomposition is well defined (see \cite{BB05}).

We define the \emph{type $B/C$ Hecke algebra} to be the associative algebra $\mathcal{H}=\mathcal{H}(B_n)$ with a linear basis $\{T_{\omega}\}_{\omega \in B_n}$ and multiplication which satisfy the following rules for any Coxeter generators $s_k$ and $\omega \in B_n$:
\begin{align}
     \label{eq: R_k left} T_{s_k} T_{\omega}&=\begin{cases}T_{s_k\omega}&\text{ if }l(s_k \omega)=l(\omega)+1,
    \\T_{\omega}&\text{ if }l(s_k\omega)=l(\omega)-1.
    \end{cases}
\end{align}
We refer to [Hum90] for the general definition of the type $B/C$ Hecke algebra, and to
[Buf20, HS23] for probabilistic applications in interacting particle systems.

Now we consider the following \emph{random walk on the type $B/C$ Hecke algebra $\mathcal{H}(B_n)$}. We place rate 1 exponential clocks at each edge $(k,k+1)$ (associated with $s_k$, $k \geq 1$), and a rate $\alpha$ exponential clock at the edge $(-1,1)$ (associated with $s_0$). For any time $t \geq 0$, we define the stochastic process $W(t)$ which takes values in $\mathcal{H}(B_n)$. Its initial value is $W(0)=T_{id}$ (the basis vector corresponding to the identity permutation). When the clock associated with $s_k$ rings at a certain time $\tau \in \mathbb{R}_{\geq 0}$, then we set $W(\tau)=T_{s_k}W(\tau^{-})$. 

Let us make the link here between the random walk $W(t)$ and the colored half-space TASEP as constructed in Section \ref{ssec:Chs-TASEP}. Note that the process $\{W(t), t \geq 0\}$, which takes values in the Hecke algebra $\mathcal{H}(B_n)$, immediately induces a stochastic process $\{\omega_{t}, t \geq 0\}$ on $B_n$. The definition of $W(t)$ and the multiplication rules \eqref{eq: R_k left} imply that the process $\{\omega_{t}, t \geq 0\}$ on $B_n$ is exactly the colored half-space TASEP as constructed in Section \ref{ssec:Chs-TASEP}: the first rule of \eqref{eq: R_k left} says that if $\omega(k)<\omega(k+1)$ and the exponential clock associated with $s_k$ rings, $\omega$ gets updated as $s_k \omega$; the second rule of \eqref{eq: R_k left} says that if $\omega(k)>\omega(k+1)$, $\omega$ remains unchanged. Note that by definition $W(0)=T_{id}$, which means $\omega_0=id$.

We record explicitly why the above dynamics may be viewed as a process on the hyperoctahedral group. 
For $n\in\mathbb N$, let
\[
\Omega_n=\{\eta:\{\pm1,\ldots,\pm n\}\to \{\pm1,\ldots,\pm n\}: 
\eta \text{ is a bijection and } \eta(-x)=-\eta(x)\}.
\]
We identify $\Omega_n$ with the type $B/C$ Weyl group $B_n$ by writing $\eta(x)$ for the color at position $x$.

For $k\geq1$, let $U_k$ be the local update which, if $\eta(k)<\eta(k+1)$, swaps the colors at $k$ and $k+1$ and simultaneously swaps the colors at $-k$ and $-k-1$; otherwise it leaves the configuration unchanged. 
Let $U_0$ be the boundary update which, if $\eta(-1)<\eta(1)$, swaps the colors at $-1$ and $1$; otherwise it leaves the configuration unchanged.

\begin{lemma}
\label{lem:preserve-hyperoctahedral}
For every $k\geq0$, the map $U_k$ sends $\Omega_n$ to itself. 
Consequently, the colored half-space TASEP preserves the hyperoctahedral structure.
\end{lemma}

\begin{proof}
It is enough to check the defining relation $\eta(-x)=-\eta(x)$. 
For $k\geq1$, the update $U_k$ exchanges the two positive positions $k,k+1$ and, at the same time, exchanges the reflected negative positions $-k,-k-1$. 
Thus, if after the update the color at $k$ is the old color at $k+1$, then the color at $-k$ is the old color at $-k-1$, which is the negative of the old color at $k+1$. 
The same verification applies to $k+1$ and $-k-1$, while all other positions are unchanged. 
For $k=0$, the update swaps the pair $-1,1$, and the relation is again preserved. 
Since each update is a composition of transpositions, bijectivity is also preserved. 
Therefore $U_k\eta\in\Omega_n$ whenever $\eta\in\Omega_n$.
\end{proof}

The following lemma makes precise the relation between the graphical construction and the Hecke-algebra product. 
We use the convention that, if the clock-ring sequence is
\[
a_1,a_2,\ldots,a_r \in \{0,1,\ldots,n-1\},
\]
then $a_1$ is the first clock to ring.

\begin{lemma}
\label{lem:clock-sequence-hecke-product}
Fix a realization of the clocks up to time $t$ for which only the clocks indexed by
$a_1,a_2,\ldots,a_r$ ring. 
Then the colored half-space TASEP configuration at time $t$, started from $\eta_0\in\Omega_n$, is
\[
U_{a_r}\cdots U_{a_2}U_{a_1}\eta_0.
\]
Equivalently, under the identification of configurations with basis elements of the type $B/C$ Hecke algebra, this configuration is represented by
\[
T_{s_{a_r}}\cdots T_{s_{a_2}}T_{s_{a_1}}T_{\eta_0}.
\]
In particular, for the packed initial condition $\eta_0=id$, the random walk $W(t)$ satisfies
\[
W(t)=T_{s_{a_r}}\cdots T_{s_{a_2}}T_{s_{a_1}}.
\]
\end{lemma}

\begin{proof}
The statement follows by induction on the number of clock rings. 
For one ring, the multiplication rule \eqref{eq: R_k left} says exactly that the update by $s_k$ is performed if the corresponding local order is increasing, and otherwise the configuration is unchanged. 
This is precisely the local rule defining $U_k$. 
Assume the claim holds after the first $r-1$ rings. 
At the $r$-th ring, the graphical construction applies $U_{a_r}$ to the current configuration, while the Hecke process left-multiplies the current basis element by $T_{s_{a_r}}$. 
Again \eqref{eq: R_k left} gives exactly the same update rule. 
This proves the claim.
\end{proof}

\begin{remark}
While the constructions in this section describe stochastic processes on a finite space, they can be extended to the infinite setting using the standard argument of Harris (\cite{Har78}). Namely, at any fixed time, the processes on an infinite space are, with probability 1, (infinite) collections of finite-space processes that did not interact with each other. Consequently, the results for processes on an infinite space follow directly from the corresponding results for their finite-space counterparts.
\end{remark}

\begin{remark}
Since the random element of the Hecke algebra and the corresponding random permutation that it induces are equivalent, for convenience, we will use the two concepts equivalently in the following descriptions and do not mention it again.
\end{remark}

The reason we describe the evolution of colored half-space TASEP as a random walk on the type $B/C$ Hecke algebra is that we wish to use a key property of the Hecke algebra $\mathcal{H}(B_n)$: the anti-involution $\iota$, which sends $\iota(T_w)=T_{w^{-1}}$ for all $w\in B_n$. In more detail, for any $T_{\omega_1},\cdots,T_{\omega_k}\in\mathcal{H}(B_n)$ we have 
\begin{align}
\label{eq:involution}
\iota(T_{\omega_n}T_{\omega_{n-1}} \cdots T_{\omega_2}T_{\omega_1})=\iota(T_{\omega_1})\iota(T_{\omega_2})\cdots \iota(T_{\omega_{n-1}})\iota(T_{\omega_n}),
\end{align}
which can be straightforwardly proved by induction in $l(\omega)$ with the use of  \eqref{eq: R_k left}. \eqref{eq:involution} was first introduced in probabilistic setting by \cite{Buf20}. Furthermore, we use the fact for any Coxeter generator $\iota(T_{s_{k}})=T_{s_{k}}$ to induce the following lemma:
\begin{lemma}
\label{lem:inverse}
For any $T_{s_{k_1}},\cdots,T_{s_{k_n}}\in\mathcal{H}(B_n)$ we have 
\begin{align}
\label{eq:involution-generator}
\iota(T_{s_{k_n}}T_{s_{k_{n-1}}} \cdots T_{s_{k_2}}T_{s_{k_1}})=T_{s_{k_1}}T_{s_{k_2}}\cdots T_{s_{k_{n-1}}}T_{s_{k_n}}.
\end{align}
\end{lemma}

\begin{remark}
In the probabilistic setting of type A Hecke algebra, Lemma \ref{lem:inverse} was proved in \cite[Lemma 2.1]{AHR09}, see also \cite[Theorem 2.1]{BF22}. The ASEP generalizations of \cite[Lemma 2.1]{AHR09} can be found in \cite[Lemma 3.1]{AAV11} and \cite[Theorem 2.2]{BB21}.
\end{remark}

By the time-reflection symmetry of a homogeneous Poisson process, we have the following lemma as the corollary of Lemma \ref{lem:inverse}:
\begin{lemma}
\label{lem:d-equality}
With $\overset{\text{d}}{=}$ denoting equality in distribution, we have 
\begin{align}
\iota(W(t)) \overset{\text{d}}{=} W(t).
\end{align}
\end{lemma}

\begin{remark}
See \cite[Lemma 3.4]{HS23} for the ASEP generalizations of Lemma \ref{lem:d-equality}. In the setting of type A Hecke algebra, Lemma \ref{lem:d-equality} was also used in \cite{AAV11,BB21,BF22,BN22}.
\end{remark}

\subsection{Color--position symmetry theorem}
\label{ssec:Symmetry-theorem}

Given a sequence of Coxeter generators
$
s_{k_1},s_{k_2},\ldots,s_{k_n},
$
we consider two generated processes. 
The order convention is that $s_{k_1}$ is applied first.
Let
$
\Theta=T_{s_{k_n}}\cdots T_{s_{k_2}}T_{s_{k_1}}.
$
This is the deterministic preprocessing word obtained by applying the updates
$s_{k_1},s_{k_2},\ldots,s_{k_n}$ to the identity configuration.

The first generated process $\eta_t^{\mathrm{gen};1}$ is obtained by first applying this deterministic word and then running the random colored half-space TASEP for time $t$. 
Thus, conditional on a clock-ring word
\[
W(t)=T_{s_{a_r}}\cdots T_{s_{a_1}},
\]
we have
\[
\eta_t^{\mathrm{gen};1}
=
T_{s_{a_r}}\cdots T_{s_{a_1}}
T_{s_{k_n}}\cdots T_{s_{k_1}}.
\]

The second generated process $\eta_t^{\mathrm{gen};2}$ is obtained in the opposite order: first run the random colored half-space TASEP for time $t$, and then apply the deterministic updates in the reverse order
\[
s_{k_n},s_{k_{n-1}},\ldots,s_{k_1}.
\]
Equivalently, with the above convention,
\[
\eta_t^{\mathrm{gen};2}
=
T_{s_{k_1}}\cdots T_{s_{k_n}}
T_{s_{a_r}}\cdots T_{s_{a_1}}.
\]
Thus $\eta_t^{\mathrm{gen};1}$ has deterministic preprocessing before the random evolution, while $\eta_t^{\mathrm{gen};2}$ has deterministic postprocessing after the random evolution. 
In the applications below, the deterministic words are chosen to be sorting words on finite intervals; hence the postprocessing description means that, after the random evolution, the colors in the specified finite blocks are sorted.

By abuse of notation, we use the same symbol $\iota$ for the induced inverse map on configurations: if $\eta$ is a signed permutation, then $\iota(\eta)=\eta^{-1}$, so that $\iota(\eta)(x)$ is the position of the particle with color $x$.

\begin{theorem}[color--position symmetry]
\label{thm:color-position}
The random configurations $\eta_t^{\mathrm{gen};1}$ and $\iota(\eta_t^{\mathrm{gen};2})$ have the same distribution. Equivalently, $\iota(\eta_t^{\mathrm{gen};1})$ and $\eta_t^{\mathrm{gen};2}$ have the same distribution.
\end{theorem}

\begin{proof}
We prove the first distributional identity; the second follows by applying the same argument after taking inverses.

Condition on a realization of the clocks up to time $t$, and write the corresponding Hecke word as
\[
W(t)=T_{s_{a_r}}\cdots T_{s_{a_1}}.
\]
By Lemma \ref{lem:clock-sequence-hecke-product} and by the definitions of the two generated processes,
\[
\eta_t^{\mathrm{gen};1}
=
T_{s_{a_r}}\cdots T_{s_{a_1}}
T_{s_{k_n}}\cdots T_{s_{k_1}},
\]
whereas
\[
\eta_t^{\mathrm{gen};2}
=
T_{s_{k_1}}\cdots T_{s_{k_n}}
T_{s_{a_r}}\cdots T_{s_{a_1}}.
\]
Applying the anti-involution $\iota$ and using $\iota(T_{s_j})=T_{s_j}$ for every Coxeter generator gives
\[
\iota(\eta_t^{\mathrm{gen};2})
=
\iota\!\left(
T_{s_{k_1}}\cdots T_{s_{k_n}}
T_{s_{a_r}}\cdots T_{s_{a_1}}
\right)
=
\iota(W(t))\,T_{s_{k_n}}\cdots T_{s_{k_1}}.
\]
By Lemma \ref{lem:d-equality}, $\iota(W(t))\overset{d}{=}W(t)$. 
Therefore
\[
\iota(\eta_t^{\mathrm{gen};2})
\overset{d}{=}
W(t)T_{s_{k_n}}\cdots T_{s_{k_1}}
=
\eta_t^{\mathrm{gen};1}.
\]
This proves
\[
\eta_t^{\mathrm{gen};1}\overset{d}{=}\iota(\eta_t^{\mathrm{gen};2}).
\]

Finally, since $\iota$ is an involution, applying $\iota$ to both sides gives
\[
\iota(\eta_t^{\mathrm{gen};1})
\overset{d}{=}
\eta_t^{\mathrm{gen};2}.
\]
\end{proof}

\begin{remark}[Intuition for the half-space color--position symmetry]
\label{rmk:intuition-color-position}
A fully colored TASEP configuration can be read in two equivalent ways. 
One may look site by site and ask which color occupies a given position, or one may look color by color and ask where a given color is located. 
Color--position symmetry says that these two descriptions have the same distribution after the order of the updates is reversed.

In the full-space multi-color TASEP, the local updates are adjacent swaps. 
Thus a sequence of clock rings produces a sequence of local comparisons and swaps. 
If, instead of tracking the color at each site, one tracks the position of each color, the same swaps are seen in the reverse order. 
This is the basic reason why a deterministic rearrangement of colors before the random evolution can be converted into a deterministic rearrangement after the random evolution, once colors and positions are interchanged.

In the half-space model, the positive and negative sides are tied together by the signed relation
\[
\eta(-x)=-\eta(x).
\]
Therefore one should not think of the two half-lines as independent copies of the full-space process. 
Knowing the colors on the positive half-line determines the colors on the negative half-line. 
The boundary update at $(-1,1)$ is the extra move that exchanges the two signed partners across the origin. 
Thus the same color--position principle still holds, but colors and positions must be exchanged as signed objects. 
This is the content of Theorem \ref{thm:color-position}.
\end{remark}

\begin{remark}
In the probabilistic setting of type A Hecke algebra, Theorem \ref{thm:color-position} was proved in \cite[Theorem 1.4]{AAV11}, \cite[Theorem 3.1]{BB21} and \cite[Proposition 2.3]{BF22}.
\end{remark}

\begin{lemma}[Generator sequences and projections]
\label{lem:generator-projection}
Let $\Pi$ be a projection of the color set into finitely many ordered classes, obtained by grouping colors into disjoint intervals. 
For a fully colored configuration $\eta$, write $\Pi\eta$ for the projected configuration. 
Let $U_k$ denote the update of the fully colored half-space TASEP corresponding to the Coxeter generator $s_k$, and let $\overline U_k$ denote the induced update on the projected multi-species configuration. 
Then, for every $k\geq 0$,
\[
\Pi(U_k\eta)=\overline U_k(\Pi\eta).
\]
Consequently, for any finite generator sequence $s_{k_1},\ldots,s_{k_m}$,
\[
\Pi(U_{k_m}\cdots U_{k_1}\eta)
=
\overline U_{k_m}\cdots \overline U_{k_1}(\Pi\eta).
\]
The same pathwise identity holds for any realization of the clock-ring sequence in the graphical construction.
\end{lemma}

\begin{proof}
The update $U_k$ depends only on the relative order of the two colors involved in the edge corresponding to $s_k$. 
Since $\Pi$ is obtained by grouping colors into ordered intervals, it preserves the order between different projected classes. 
If the two colors belong to different projected classes, then the fully colored update swaps them if and only if the corresponding projected update swaps the two classes. 
If the two colors belong to the same projected class, the fully colored update may swap the two colors, but this swap is invisible after applying $\Pi$. 
Thus $\Pi(U_k\eta)=\overline U_k(\Pi\eta)$ for each $k\geq1$.

The same argument applies to the boundary generator $s_0$, whose update compares the colors at $-1$ and $1$. 
Therefore the identity holds for all Coxeter generators. 
Iterating the one-step identity along the sequence $s_{k_1},\ldots,s_{k_m}$ gives the formula for deterministic generator sequences. 
Applying the same argument to the random sequence of clock rings gives the pathwise identity for the graphical construction.
\end{proof}

\section{Exact distribution of the second-class particles in the half-line open TASEP with shocks}
\label{sec:distributional-identity}
In this section, we use the color--position symmetry theorem introduced in Section \ref{sec:ColorPositionSymmetry}
 to derive the exact distribution of the second-class particles in different kinds of multi-color half-line open TASEP. We focus on the one-shock case from Section \ref{ssec:one-shock} to Section \ref{ssec:Distribution-Identity}, and discuss the two-shock case in Section \ref{ssec:two-shock}.

\subsection{Second-class particles in one shock}
\label{ssec:one-shock}
Let $M_1, M_2$ be positive integers, and let $\pi_{M_1,M_2}$ be a signed permutation of the set $\{-M_1-M_2-1, \cdots, -1, 1, \cdots, M_1+M_2+1\}$ such that $\pi_{M_1,M_2}(i)=M_1+M_2+2-i$ and $\pi_{M_1,M_2}(-i)=-\pi_{M_1,M_2}(i)$ for $i=1,2, \cdots, M_1+M_2+1$. We assume the length $l(\pi_{M_1,M_2})=m$ and let $\pi_{M_1,M_2}=s_{k_m} s_{k_{m-1}} \cdots s_{k_1}$ be a minimal length decomposition of $\pi_{M_1,M_2}$ into the Coxeter generators (there can be many such decompositions, we only choose one of them).

Consider the first process $\eta_{t}^{gen;1}$ introduced in Section \ref{ssec:Symmetry-theorem} with generator sequence $s_{k_1}, s_{k_2}, \cdots,s_{k_m}$. In this process, we start from the configuration $T_{id}$ and sequentially apply updates $T_{s_{k_1}}, T_{s_{k_2}}, \cdots, T_{s_{k_m}}$ to it. This can be expressed as $\eta_{0}^{gen;1}=T_{s_{k_m}} T_{s_{k_{m-1}}}\cdots T_{s_{k_1}} T_{id}$. To elaborate, within the packed initial configuration, we individually arrange integers within the ranges $[1, \cdots, M_1+M_2+1]$ and $[-M_1-M_2-1, \cdots, -1]$ in reverse order. Subsequently, we run the continuous-time dynamics to obtain the random configuration $\eta_{t}^{gen;1}$ up to time $t$.

In this section, we will consider two kinds of tricolor half-line open TASEP which can be regarded as different projections of $\eta_{t}^{gen;1}$. 

Firstly, we consider a tricolor half-line open TASEP denoted as $\eta_{t}$ with the initial condition
\begin{align}
 \label{eq:eta_0}
\eta_{0}(x)&=
\begin{cases}
1, & \text{ if } x \leq -1 \ \text{or}\  M_1+2 \leq x \leq M_1+M_2+1, \\
2, & \text{ if } x=M_1+1, \\
+\infty, & \text{otherwise}.
\end{cases}
\end{align}
The interpretation of $\eta_{t}$ in terms of $\eta_{t}^{gen;1}$ (with Coxeter-generator sequence $s_{k_1}, s_{k_2}, \cdots,s_{k_m}$) is the following. Since $\pi_{M_1,M_2}(M_1+1)=M_2+1$, we identify the particle with color $M_2+1$ as the second-class particle (denoted by 2). Furthermore, particles with colors less than $M_2+1$ are called first-class particles (denoted by 1), and particles with colors greater than $M_2+1$ are holes (denoted by $+\infty$). Secondly, we consider a tricolor half-line open TASEP denoted as $\tilde{\eta}_{t}$ with the following initial condition
\begin{align}
 \label{eq:tilde-eta_0}
\tilde{\eta}_{0}(x)&=
\begin{cases}
1, & \text{ if } x \leq -1,  \\
2, & \text{ if } M_1+1 \leq x \leq M_1+M_2+1, \\
+\infty, & \text{otherwise}.
\end{cases}
\end{align}
Similar to $\eta_{t}$, we can interpret $\tilde{\eta}_{t}$ in terms of $\eta_{t}^{gen;1}$ (with Coxeter-generator sequence $s_{k_1}, s_{k_2}, \cdots,s_{k_m}$) as follows: we identify the particles with negative colors as the first-class particles, the particles with colors in $[1, M_2+1]$ as second-class particles, and the particles with colors greater than $M_2+1$ as holes. 

\begin{remark}
\label{rem:projection-sorting}
By Lemma \ref{lem:generator-projection}, each of the projections used in Section \ref{sec:distributional-identity} commutes with the Coxeter-generator updates. 
Consequently, if a fully colored process is driven by a deterministic or random generator sequence, then its projection agrees pathwise with the corresponding multi-species half-line open TASEP driven by the same generator sequence. 
Thus Theorem \ref{thm:color-position} may be applied first at the fully colored level and then projected to the multi-species processes used in this section.
\end{remark}

\begin{remark}
\label{rmk:model-equi}
The half-space dynamics preserve the signed symmetry between $\mathbb N$ and $-\mathbb N$. 
Since the negative half-line is filled with first-class particles in the initial configurations \eqref{eq:eta_0} and \eqref{eq:tilde-eta_0}, the restriction of the half-space processes $\eta_t$ and $\tilde{\eta}_t$ to the positive half-line can be identified with the corresponding half-line open TASEP processes shown in the upper and lower panels of Figure \ref{fig:one-shock}, respectively. 
For this reason, we keep the same notation. 
The same convention is used for the two-shock case in Section \ref{ssec:two-shock}.
\end{remark}

\subsection{Standard half-line open TASEP with step initial condition}
\label{ssec:HL-TASEP}

\begin{figure}[h]
\begin{tikzpicture}[scale=1.25]
\node[circle, inner sep=10pt, draw=black, fill=orange, opacity=0.7] (a) at (-3.3, 0) {};
\draw (-3.3, 0) node[scale=0.6,blue] {$\text{reservoir}$};
\node[circle, inner sep=3.8pt, draw=black, fill=white, opacity=0.7] (b) at (-2, 0) {};
\node[circle, inner sep=3.8pt, draw=black, fill=white, opacity=0.7] (c) at (-1, 0) {};
\node[circle, inner sep=3.8pt, draw=black, fill=white, opacity=0.7] (d) at (0, 0) {};
\node[circle, inner sep=3.8pt, draw=black, fill=white, opacity=0.7] (e) at (1, 0) {};
\node[circle, inner sep=3.8pt, draw=black, fill=white, opacity=0.7] (f) at (2, 0) {};
%
\draw (a) -- (b);
\draw (b) -- (c);
\draw (c) -- (d);
\draw (d) -- (e);
\draw (e) -- (f);
\draw[dotted, thick] (f) -- (3,0);
\draw (-2, -0.5) node[scale=0.5] {$1$};
\draw (-1, -0.5) node[scale=0.5] {$2$};
\draw (0, -0.5) node[scale=0.5] {$3$};
\draw (1, -0.5) node[scale=0.5] {$4$};
\draw (2, -0.5) node[scale=0.5] {$\cdots$};
\draw[-stealth, dashed, semithick] (-3.3,0.4) to[out=45,in=135] node[midway, above] {$\alpha$} (-2,0.2) ;
\end{tikzpicture}
\caption{The standard half-line open TASEP with step initial condition.}\label{fig:standard}
\end{figure}
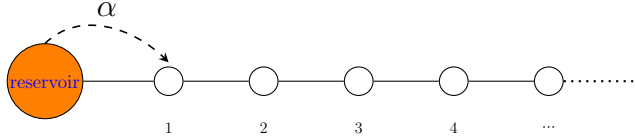

In this section, we introduce a two-color half-line open TASEP $\eta_{t}^{step}(x)$ starting from the step initial condition
\begin{align}
 \label{eq:eta-step}
\eta_{0}^{step}(x)&=
\begin{cases}
1, & \text{ if } x \leq -1, \\
+\infty, & \text{ if } x \geq 1.
\end{cases}
\end{align}
We can interpret $\eta_{t}^{step}(x)$ as the projection of the colored half-space TASEP introduced in Section \ref{ssec:Chs-TASEP}. In this context, if a color is $\leq -1$, it signifies the presence of a particle (denoted by 1), whereas a color $\geq 1$ indicates the existence of a hole (denoted by $+\infty$). One can identify the two-color half-line open TASEP $\eta_{t}^{step}(x)$ with the standard half-line open TASEP in Figure \ref{fig:standard}. For any $x \in \mathbb{R}_{\geq 0}$, we use $\mathcal{N}(x,t)$ to represent the number of particles located weakly to the right of $x$ in $\eta_{t}^{step}$, namely,
\begin{align}
\mathcal{N}(x,t)=\sum_{z \geq x} \delta_{\eta_t^{step}(z), 1}.
\end{align}

\subsection{Distribution identity}
\label{ssec:Distribution-Identity}
In this section, we study the quantities $f(t)$, $\mathbf{N}_{1}(x,t)$ and $\mathbf{N}_{2}(x,t)$ by relating them to $\mathcal{N}(x,t)$ via the color--position symmetry theorem \ref{thm:color-position}. 
\begin{theorem}
\label{thm:tricolor-disappear}
The probability of the unique second-class particle disappearing from the half-line open TASEP $\eta_t$ is given by:
\begin{align}
\label{eq:tricolor-disappear}
\mathbb{P}(f(t)<0)=\mathbb{P}(\mathcal{N}(1,t)-\mathcal{N}(M_1+M_2+2,t) \geq M_1+1).
\end{align}
Equivalently, the probability that the second-class particle exists in the system is equal to
\begin{align}
\label{eq:tricolor-exist}
\mathbb{P}(f(t)>0)=\mathbb{P}(\mathcal{N}(1,t)-\mathcal{N}(M_1+M_2+2,t) \leq M_1).
\end{align}
\end{theorem}

\begin{proof}
We use Coxeter-generator sequence $s_{k_1}, s_{k_2}, \cdots,s_{k_m}$ introduced in Section \ref{ssec:one-shock} in the following. By the relation of $\eta_t$ and $\eta_{t}^{gen;1}$, we have
\begin{align*}
\mathbb{P}(f(t)<0)=\mathbb{P}(\iota( \eta_{t}^{gen;1})(M_2+1)<0);
\end{align*}
Using Theorem \ref{thm:color-position}, we have
\begin{align*}
\mathbb{P}(\iota( \eta_{t}^{gen;1})(M_2+1)<0)=\mathbb{P}(\eta_{t}^{gen;2}(M_2+1)<0),
\end{align*}
where $\eta_{t}^{gen;2}(M_2+1)<0$ means the color at position $M_2+1$ is less than $0$ in the random configuration $\eta_{t}^{gen;2}$. Recall that we get the random configuration $\eta_{t}^{gen;2}$ as follows: we start from the identity configuration, and then we run the continuous-time process for time $t$, and finally we sort the colors on the interval $[1,M_1+M_2+1]$ in a decreasing order, at the same time we sort the colors on the interval $[-M_1-M_2-1,-1]$ in a decreasing order.

Consider the projection introduced in Section \ref{ssec:HL-TASEP}: colors $\leq -1$ are interpreted as particles, while colors $\geq 1$ are interpreted as holes. 
Under this projection, before the final deterministic sorting step, the process appearing in $\eta_t^{\mathrm{gen};2}$ has the same law as the standard half-line open TASEP with step initial condition.

Let
\[
I=\{1,\ldots,M_1+M_2+1\}, \qquad L=|I|=M_1+M_2+1,
\]
and let
\[
P_I(t)=\mathcal N(1,t)-\mathcal N(M_1+M_2+2,t)
\]
be the number of particles in $I$ at time $t$ for the standard half-line open TASEP. 
The postprocessing in $\eta_t^{\mathrm{gen};2}$ sorts the colors on $I$ in decreasing order. 
After the particle/hole projection, this means that the holes in $I$ are placed on the left and the particles in $I$ are placed on the right. 
Hence, after sorting, the first $L-P_I(t)$ sites of $I$ are holes and the last $P_I(t)$ sites are particles.

Therefore the event $\eta_t^{\mathrm{gen};2}(M_2+1)<0$ occurs if and only if the site $M_2+1$ belongs to the particle block after sorting (see Figure \ref{fig:one-shock-sorting} for illustration). 
This is equivalent to
\[
M_2+1>L-P_I(t),
\]
or, since $L=M_1+M_2+1$, to
\[
P_I(t)\geq M_1+1.
\]
Substituting the definition of $P_I(t)$ gives \eqref{eq:tricolor-disappear}.
\end{proof}

\begin{theorem}
\label{thm:tricolor-height}
We have the following distribution identity relating $\mathbf{N}_{1}(1,t)$, $\mathbf{N}_{2}(1,t)$ and $\mathcal{N}(x,t)$:
\begin{multline}
\label{eq:remain}
\left(\mathbf{N}_{1}(1,t) , \mathbf{N}_{2}(1,t) \right) \stackrel{d}{=} \\ 
\left( \mathcal{N}(1, t) , \min \left\{M_1+M_2+1- \mathcal{N}(1, t) + \mathcal{N}(M_1 + M_{2}+2, t), M_{2}+1 \right\} \right),
\end{multline}
where $\stackrel{d}{=}$ denotes equality in distribution.
\end{theorem}

\begin{proof}
We use the same projection as in the proof of Theorem \ref{thm:tricolor-disappear}. 
By the relation of $\tilde\eta_t$ and $\eta_{t}^{gen;1}$, we have
\begin{multline*}
\left(\mathbf{N}_{1} (1,t) , \mathbf{N}_{2} (1,t) \right) \stackrel{d}{=} \\ 
\left( \#\{i: i<0, \iota( \eta_{t}^{gen;1})(i) \geq 1\},  \#\{i: 1 \leq i \leq M_2+1, \iota( \eta_{t}^{gen;1})(i) \geq 1\}\right),
\end{multline*}
Using Theorem \ref{thm:color-position}, we have
\begin{multline*}
\left( \#\{i: i<0, \iota( \eta_{t}^{gen;1})(i) \geq 1\},  \#\{i: 1 \leq i \leq M_2+1, \iota( \eta_{t}^{gen;1})(i) \geq 1\}\right) \stackrel{d}{=} \\
\left( \#\{i: i<0, \eta_{t}^{gen;2}(i) \geq 1\},  \#\{i: 1 \leq i \leq M_2+1, \eta_{t}^{gen;2}(i) \geq 1\}\right),
\end{multline*}
We now translate the two coordinates on the right-hand side into particle counts for the standard half-line open TASEP. 
For the first coordinate, by the signed symmetry of the fully colored half-space process, the number of negative colors whose positions are positive equals the number of particles on the positive half-line under the particle/hole projection. 
Therefore
\[
\#\{i: i<0,\eta_t^{\mathrm{gen};2}(i)\geq 1\}
=
\mathcal N(1,t).
\]

For the second coordinate, set again
\[
I=\{1,\ldots,M_1+M_2+1\},\qquad L=M_1+M_2+1,
\]
and
\[
P_I(t)=\mathcal N(1,t)-\mathcal N(M_1+M_2+2,t).
\]
Thus $P_I(t)$ is the number of particles in $I$, and
\[
H_I(t)=L-P_I(t)
=
M_1+M_2+1-\mathcal N(1,t)+\mathcal N(M_1+M_2+2,t)
\]
is the number of holes in $I$. 
In $\eta_t^{\mathrm{gen};2}$, the postprocessing sorts the interval $I$ so that, after the particle/hole projection, all holes in $I$ are placed to the left and all particles in $I$ are placed to the right. 
Consequently, among the first $M_2+1$ sites of $I$, the number of holes is
$\min\{H_I(t),M_2+1\}$ (see Figure \ref{fig:one-shock-sorting} for illustration).
Equivalently,
\[
\#\{1\leq i\leq M_2+1:\eta_t^{\mathrm{gen};2}(i)\geq 1\}
=
\min \left\{
M_1+M_2+1-\mathcal N(1,t)+\mathcal N(M_1+M_2+2,t),
M_2+1
\right\}.
\]
Together with the color--position symmetry identity above, this proves \eqref{eq:remain}.
\end{proof}

\begin{figure}[t]
\centering
\begin{tikzpicture}[scale=0.9]
    \draw[thick] (0,0) -- (12,0);

    \draw[thick] (0,-0.15) -- (0,0.15);
    \draw[thick] (12,-0.15) -- (12,0.15);

    \node[below] at (0,-0.15) {$1$};
    \node[below] at (12.5,-0.15) {$M_1+M_2+1$};

    \draw[decorate,decoration={brace,amplitude=5pt}] (0,0.45) -- (5,0.45);
    \node[above] at (2.5,0.65) {$H_I(t)=L-P_I(t)$ holes};

    \draw[decorate,decoration={brace,amplitude=5pt}] (5,0.45) -- (12,0.45);
    \node[above] at (8.5,0.65) {$P_I(t)$ particles};

    \foreach \x in {0.5,1.2,1.9,2.6,3.3,4.0,4.7}
        \draw[fill=white] (\x,0) circle (3pt);

    \foreach \x in {5.4,6.1,6.8,7.5,8.2,8.9,9.6,10.3,11.0,11.7}
        \draw[fill=black] (\x,0) circle (3pt);

    \draw[dashed,thick] (5.0,-0.75) -- (5.0,1.15);
    \node[below] at (5.0,-0.75) {$M_2+1$};

    \node at (2.5,-0.55) {holes after sorting};
    \node at (8.5,-0.55) {particles after sorting};

\end{tikzpicture}
\caption{Sorting of colors in the proof of Theorems \ref{thm:tricolor-disappear} and \ref{thm:tricolor-height}. 
Under the particle/hole projection, the postprocessing in $\eta_t^{\mathrm{gen};2}$ sorts the interval $I=\{1,\ldots,M_1+M_2+1\}$ by placing holes to the left and particles to the right. 
The site $M_2+1$ lies in the particle block if and only if the number $P_I(t)$ of particles in $I$ is at least $M_1+1$. 
The number of holes among the first $M_2+1$ sites is $\min\{H_I(t),M_2+1\}$.}
\label{fig:one-shock-sorting}
\end{figure}

\begin{example}
Let $M_1=1$ and $M_2=2$. 
Then the relevant sorting interval in Theorem \ref{thm:tricolor-disappear} is
\[
I=\{1,2,3,4\}.
\]
The quantity
\[
\mathcal N(1,t)-\mathcal N(5,t)
\]
is the number of particles in $I$ at time $t$ for the standard half-line open TASEP. 
After the deterministic sorting step in $\eta_t^{\mathrm{gen};2}$, holes in $I$ are placed to the left and particles to the right. 
The marked site is $M_2+1=3$. 
Hence this site is occupied by a particle exactly when at least two particles are present in $I$, namely when
\[
\mathcal N(1,t)-\mathcal N(5,t)\geq 2.
\]
This is the special case $M_1=1$, $M_2=2$ of \eqref{eq:tricolor-disappear}.
\end{example}

\begin{remark}
Theorem \ref{thm:tricolor-disappear} and Theorem \ref{thm:tricolor-height} are similar to the results in \cite[Proposition 6.1]{BB21} and \cite[Proposition 2.5]{BF22}, but in the half-space geometry. Unlike in the full-space geometry, we lose the property of translation-invariance. That is why \cite[Proposition 2.5]{BF22} can give an accurate description of the distribution function, while we only obtain the probability distribution around the origin.
\end{remark}

\begin{remark}
We can generalize the results in Theorem \ref{thm:tricolor-disappear} and Theorem \ref{thm:tricolor-height} to the distribution identities of distribution functions by introducing the half-line open TASEP with the shifted step initial condition (denoted by $\eta_{t}^{\mathrm{step}(y)}$, $y \in \mathbb{N} \cup -\mathbb{N}$): 
\begin{align*}
\eta_{0}^{\mathrm{step}(y)}(x)&=
\begin{cases}
1, & \text{ if } x \leq y, \\
+\infty, & \text{ if } x > y.
\end{cases}
\end{align*}
where $x \in \mathbb{N} \cup -\mathbb{N}$. We can interpret $\eta_{t}^{\mathrm{step}(y)}$ as the projection of the colored half-space TASEP introduced in Section \ref{ssec:Chs-TASEP}: if a color is $\leq y$, we identify it as a particle (denoted by 1), whereas a color $>y$ is regarded as a hole (denoted by $+\infty$). For any $x \in \mathbb{R}$, we use $\mathcal{N}^{y}(x,t)$ to represent the number of particles which are weakly to the right of $x$ in $\eta_{t}^{\mathrm{step}(y)}$, namely,
\begin{align*}
\mathcal{N}^{y}(x,t)=\sum_{z \geq x} \delta_{\eta_t^{\mathrm{step}(y)}(z), 1}.
\end{align*}
Then the equation \eqref{eq:tricolor-disappear} can be generalized to the following form:
\begin{align*}
\mathbb{P}(f(t) \leq y)=\mathbb{P}(\mathcal{N}^{y}(1,t)-\mathcal{N}^{y}(M_1+M_2+2,t)\geq M_1+1).
\end{align*}
The equation \eqref{eq:remain} can also be generalized to the following form:
\begin{multline*}
\left(\mathbf{N}_{1}(y,t) , \mathbf{N}_{2}(y,t) \right) \stackrel{d}{=} \\ 
\left( \mathcal{N}^{y}(1, t) -y, \min \left\{M_1+M_2+1- \mathcal{N}^{y}(1, t) + \mathcal{N}^{y} (M_1 + M_{2}+2, t), M_{2}+1 \right\} \right).
\end{multline*}
We will address the asymptotic analysis in our future work.
\end{remark}

\begin{remark}
One can also consider generalizing the results in Theorem \ref{thm:tricolor-disappear} and Theorem \ref{thm:tricolor-height} to the distribution identities of distribution functions by introducing the shift-invariance of colored half-space TASEP. Let $h^{\leq i}(x,t)$ denote the number of particles of color $\leq i$ which are weakly to the right of $x$. Then we have 
\begin{align*}
\mathbb{P}(f(t) \leq y)=\mathbb{P}(h^{\leq y}(1,t)-h^{\leq y}(M_1+M_2+1,t) \geq M_1+1).
\end{align*}
We need to derive the shift-invariance theorem to set up the connection between $h^{\leq y}(x,t)$ and $h^{\leq -1}(1,t)$, for any $x,y \in \mathbb{N}\cup -\mathbb{N}$. We have noticed that there exists specific shift-invariance theorem in \cite{He22}, but it does not fit our situation.
\end{remark}

\subsection{Second-class particles in two shocks}
\label{ssec:two-shock}
In this section, we consider two more multi-color half-line open TASEPs which have two-shock initial conditions. Let $N,M$ be positive integers. The first multi-color half-line open TASEP has three colors, we denote it as $\xi_{t}$ and set the following initial condition:
\begin{align}
 \label{eq:xi}
\xi_{0}(x)&=
\begin{cases}
1, & \text{ if } x \leq -1 \ \text{or}\ \ N+1 \leq x \leq N+M, \\
2, & \text{ if } x=N+2M+1, \\
1, & \text{ if } N+2M+2 \leq x \leq 2N+2M+1, \\
+\infty, & \text{otherwise}.
\end{cases}
\end{align}
The second multi-color half-line open TASEP has four colors, we denote it as $\tilde\xi_{t}$ and assign the following initial condition:
\begin{align}
 \label{eq:tilde-xi}
\tilde\xi_{0}(x)&=
\begin{cases}
1, & \text{ if } x \leq -1,  \\
2, & \text{ if } N+1 \leq x \leq N+M, \\
3, & \text{ if } 2M+N+1 \leq x \leq 2M+2N+1, \\
+\infty, & \text{otherwise}.
\end{cases}
\end{align}
For the same reason as in Remark \ref{rmk:model-equi}, the model $\xi_{t}$ is equivalent to the half-line open TASEP shown in the upper panel of Figure \ref{fig:two-shocks}, and the model $\tilde\xi_{t}$ is equivalent to the one shown in the lower panel of Figure \ref{fig:two-shocks}. We use the same notation $g(t)$ and $\mathsf{N}_{i}(x,t)$ as defined in Section \ref{ssec:models}.

\begin{theorem}
\label{thm:two-shock}
The probability of the unique second-class particle disappearing from the half-line open TASEP $\xi_t$ is given by:
\begin{multline}
\label{eq:four-color-disappear}
\mathbb{P} \left( g(t)<0 \right) = \mathbb{P} \left( \mathcal{N}(N+M+1,t)-\mathcal{N}(2N+2M+2,t)\right.\\
\left.+\min\{N, \mathcal{N}(1,t)-\mathcal{N}(N+M+1,t)\}\ge N+M+1 \right);
\end{multline}
We also have the following distribution identity among $\mathsf{N}_{1}(1,t)$, $\mathsf{N}_{2}(1,t)$, $\mathsf{N}_{3}(1,t)$  and $\mathcal{N}(x,t)$:
\begin{multline}
\label{eq:four-color-remain}
\left(\mathsf{N}_{1}(1,t) , \mathsf{N}_{2}(1,t), \mathsf{N}_{3}(1,t)  \right) \stackrel{d}{=} \bigl( \mathcal{N}(1, t) , \min\{M, M+N-\mathcal{N}(1,t)+\mathcal{N}(M+N+1,t)\},\bigr.\\
\bigl.\min \bigl\{N+1, M+N+1-\mathcal{N}(M+N+1,t)+\mathcal{N}(2M+2N+2,t)\bigr.\bigr.\\
\bigl.\bigl.+\max\{0, N-\mathcal{N}(1,t)+\mathcal{N}(M+N+1,t)\}\bigr\}\bigr).
\end{multline}
\end{theorem}

\begin{proof}

Let $\pi_1$ be a signed permutation of the set $\{-2M-2N-1, \cdots, -M-1, M+1, \cdots, 2M+2N+1\}$ such that $\pi_1(M+1+i)=2M+2N+1-i$ and $\pi_1(-i)=-\pi_1(i)$, where $i=0,1,\cdots, M+2N$. And let $\pi_2$ be a signed permutation of the set $\{-M-N, \cdots, -1, 1, \cdots, M+N\}$ such that $\pi_2(1+i)=M+N-i$ and $\pi_2(-i)=-\pi_2(i)$, where $i=0,1,\cdots, M+N-1$. Consider a minimal length decomposition of permutation $\pi_2 \cdot \pi_1$ into Coxeter generators: $\pi_2 \cdot \pi_1=s_{\ell_p}s_{\ell_{p-1}} \cdots s_{\ell_1}$, where $l(\pi_2 \cdot \pi_1)=p$.

We use the data $s_{\ell_1}, s_{\ell_2}, \cdots,s_{\ell_p}$ to construct the first process $\eta_{t}^{gen;1}$ introduced in Section \ref{ssec:Symmetry-theorem}.  We obtain $\xi_{t}$ from $\eta_{t}^{gen;1}$ (with Coxeter-generator sequence $s_{\ell_1}, s_{\ell_2}, \cdots,s_{\ell_p}$) by identifying the particle with color $M+N+1$ as the second-class particle (denoted by 2), particles with colors less than $M+N+1$ as first-class particles (denoted by 1), and particles with colors greater than $M+N+1$ as holes (denoted by $+\infty$). This implies:
\begin{align*}
\mathbb{P}(g(t)<0)=\mathbb{P}(\iota( \eta_{t}^{gen;1})(M+N+1)<0).
\end{align*}
Additionally, $\tilde\xi_{t}$ can be interpreted as the projection of $\eta_{t}^{gen;1}$(with Coxeter-generator sequence $s_{\ell_1}, s_{\ell_2}, \cdots,s_{\ell_p}$) as follows: we identify the particles with negative colors as first-class particles, colors in $[1, M]$ as second-class particles, colors in $[M+1, M+N+1]$ as third-class particles (denoted by 3), and colors greater than $M+N+1$ as holes. This gives:
\begin{multline*}
\left(\mathsf{N}_{1} (1,t) , \mathsf{N}_{2} (1,t), \mathsf{N}_{3} (1,t) \right) \stackrel{d}{=} \left( \#\{i: i<0, \iota( \eta_{t}^{gen;1})(i) \geq 1\}, \right. \\ 
\left. \#\{i: 1 \leq i \leq M, \iota( \eta_{t}^{gen;1})(i) \geq 1\},  \#\{i: M+1 \leq i \leq M+N+1, \iota( \eta_{t}^{gen;1})(i) \geq 1\}\right).
\end{multline*}
Theorem \ref{thm:color-position} implies:
\begin{align*}
\mathbb{P}(g(t)<0)=\mathbb{P}(\eta_{t}^{gen;2}(M+N+1)<0),
\end{align*}
and the following relation:
\begin{multline}
\label{eq:joint}
\left(\mathsf{N}_{1} (1,t) , \mathsf{N}_{2} (1,t), \mathsf{N}_{3} (1,t) \right) \stackrel{d}{=} \left( \#\{i: i<0, \eta_{t}^{gen;2}(i) \geq 1\}, \right. \\ 
\left. \#\{i: 1 \leq i \leq M, \eta_{t}^{gen;2}(i) \geq 1\},  \#\{i: M+1 \leq i \leq M+N+1, \eta_{t}^{gen;2}(i) \geq 1\}\right),
\end{multline}
where we use Coxeter-generator sequence $s_{\ell_1}, s_{\ell_2}, \cdots,s_{\ell_p}$ to construct $\eta_{t}^{gen;2}$. Recall that we get the random configuration $\eta_{t}^{gen;2}$ as follows: we start from the identity configuration; and then we run the continuous-time process for time $t$; after that we separately sort the colors on the interval $[1,M+N]$ and $[-M-N,-1]$ in a decreasing order at the same time; finally we separately sort the colors on the interval $[M+1,2M+2N+1]$ and $[-2M-2N-1,-M-1]$ in a decreasing order.

Consider the projection introduced in Section \ref{ssec:HL-TASEP}: colors $\leq -1$ are interpreted as particles, while colors $\geq 1$ are interpreted as holes. 
Under this projection, before the deterministic postprocessing, the process has the same law as the standard half-line open TASEP with step initial condition.

The postprocessing consists of two sorting operations. 
First, the interval
\[
A=\{1,\ldots,M+N\}
\]
is sorted, placing holes to the left and particles to the right. 
Second, the interval
\[
B=\{M+1,\ldots,2M+2N+1\}
\]
is sorted in the same way. 
The overlap of these two intervals is
\[
C=A\cap B=\{M+1,\ldots,M+N\}, \qquad |C|=N.
\]

Let
\[
P_A(t)=\mathcal N(1,t)-\mathcal N(M+N+1,t)
\]
be the number of particles in $A$ before the first sorting. 
After sorting $A$, the particles in $A$ occupy the rightmost positions of $A$. 
Hence the number of particles in the overlap $C$ after the first sorting is
\[
\min\{N,P_A(t)\}.
\]
The remaining part of $B$ to the right of $C$ is
\[
R=\{M+N+1,\ldots,2M+2N+1\},
\]
and the number of particles in $R$ is
\[
P_R(t)=\mathcal N(M+N+1,t)-\mathcal N(2M+2N+2,t).
\]
Therefore, just before the second sorting, the total number of particles in $B$ is
\[
P_R(t)+\min\{N,P_A(t)\}.
\]

After sorting $B$, holes occupy the left part of $B$ and particles occupy the right part. 
The marked site $M+N+1$ is the $(N+1)$-st site of $B$, or equivalently there are $M+N+1$ sites from $M+N+1$ to the right endpoint of $B$. 
Thus $\eta_t^{\mathrm{gen};2}(M+N+1)<0$ holds if and only if the sorted interval $B$ contains at least $M+N+1$ particles. 
This is equivalent to
\[
P_R(t)+\min\{N,P_A(t)\}\geq M+N+1.
\]
Substituting the definitions of $P_A(t)$ and $P_R(t)$ gives \eqref{eq:four-color-disappear}.

We now derive the three coordinates in \eqref{eq:four-color-remain}. 
The first coordinate follows from the signed symmetry of the colored half-space process:
\[
\#\{i:i<0,\eta_t^{\mathrm{gen};2}(i)\geq1\}
=
\#\{i:i>0,\eta_t^{\mathrm{gen};2}(i)\leq -1\}
=
\mathcal N(1,t).
\]

For the second coordinate, recall that after the first sorting of
\[
A=\{1,\ldots,M+N\},
\]
holes are placed to the left and particles to the right. 
The number of particles in $A$ is
\[
P_A(t)=\mathcal N(1,t)-\mathcal N(M+N+1,t),
\]
and hence the number of holes in $A$ is
\[
H_A(t)=M+N-P_A(t)
=
M+N-\mathcal N(1,t)+\mathcal N(M+N+1,t).
\]
The interval $\{1,\ldots,M\}$ is not affected by the second sorting step. 
Thus the number of holes among its $M$ sites is
\[
\min\{M,H_A(t)\}
=
\min\{M,M+N-\mathcal N(1,t)+\mathcal N(M+N+1,t)\}.
\]
This gives the second coordinate in \eqref{eq:four-color-remain}.

It remains to compute the third coordinate. 
Let
\[
C=\{M+1,\ldots,M+N\}
\]
be the overlap of the two sorting intervals. 
After sorting $A$, the number of holes in $C$ is
\[
\max\{0,N-P_A(t)\}
=
\max\{0,N-\mathcal N(1,t)+\mathcal N(M+N+1,t)\}.
\]
Let
\[
R=\{M+N+1,\ldots,2M+2N+1\}.
\]
The number of holes in $R$ is
\[
H_R(t)
=
M+N+1-\mathcal N(M+N+1,t)+\mathcal N(2M+2N+2,t).
\]
Therefore, before the second sorting, the total number of holes in
\[
B=\{M+1,\ldots,2M+2N+1\}
\]
is
\[
H_R(t)+\max\{0,N-P_A(t)\}.
\]
After sorting $B$, holes occupy the left side of $B$. 
The interval $\{M+1,\ldots,M+N+1\}$ consists of the first $N+1$ sites of $B$. 
Hence the number of holes in this interval is
\[
\min\bigl\{N+1,
H_R(t)+\max\{0,N-P_A(t)\}
\bigr\}.
\]
Substituting the expressions for $H_R(t)$ and $P_A(t)$ gives the third coordinate in \eqref{eq:four-color-remain}. 
Combining the three coordinate identities with \eqref{eq:joint} proves \eqref{eq:four-color-remain}. See Figure \ref{fig:two-shock-sorting} for illustration.
\end{proof}

\begin{figure}[t]
\centering
\begin{tikzpicture}[scale=0.85]
    \draw[thick] (0,0) -- (14,0);

    \draw[thick] (0,-0.12) -- (0,0.12);
    \draw[thick] (6,-0.12) -- (6,0.12);
    \draw[thick] (9,-0.12) -- (9,0.12);
    \draw[thick] (14,-0.12) -- (14,0.12);

    \node[below] at (0,-0.15) {$1$};
    \node[below] at (6,-0.15) {$M+1$};
    \node[below] at (9,-0.15) {$M+N+1$};
    \node[below] at (14,-0.15) {$2M+2N+1$};

    \draw[decorate,decoration={brace,amplitude=5pt}] (0,0.6) -- (9,0.6);
    \node[above] at (4.5,0.8) {$A=\{1,\ldots,M+N\}$, sorted first};

    \draw[decorate,decoration={brace,amplitude=5pt,mirror}] (6,-0.65) -- (14,-0.65);
    \node[below] at (10,-0.9) {$B=\{M+1,\ldots,2M+2N+1\}$, sorted second};

    \draw[very thick] (6,0) -- (9,0);
    \node[above] at (7.5,0.15) {$C$};

    \draw[very thick,dashed] (9,0) -- (14,0);
    \node[above] at (11.5,0.15) {$R$};

    \draw[dashed,thick] (9,-1.2) -- (9,1.25);

\end{tikzpicture}
\caption{The two sorting intervals in the proof of Theorem \ref{thm:two-shock}. 
The first sorting acts on $A=\{1,\ldots,M+N\}$ and the second on $B=\{M+1,\ldots,2M+2N+1\}$. 
Their overlap is $C=\{M+1,\ldots,M+N\}$, of size $N$, and the right part of $B$ is $R=\{M+N+1,\ldots,2M+2N+1\}$. 
After the first sorting, the number of particles in $C$ is $\min\{N,P_A(t)\}$. 
Thus, before the second sorting, the number of particles in $B$ is $P_R(t)+\min\{N,P_A(t)\}$.}
\label{fig:two-shock-sorting}
\end{figure}

\section{Asymptotic distribution of second-class particles under constant scaling}
\label{sec:constant-scaling}
In this section, we fix the size of the shocks as constants, and consider the asymptotic  distribution of the second-class particles introduced in Section \ref{sec:distributional-identity}.  In Section \ref{ssec:stationary-measure}, we review the stationary measure  of the standard half-line open TASEP with the step initial condition, and derive a combinatorial formula for the  stationary measure at $\alpha > \frac{1}{2}$. In Section \ref{ssec:constant-formulas} we give the exact formulas of the asymptotic behavior of the second-class particles.

\subsection{Stationary measure of half-line open TASEP}
\label{ssec:stationary-measure}

Let the parameter $\alpha \in (0,1)$ represent the particle density imposed by the reservoir at site 1. When $\alpha \leq \frac{1}{2}$, Liggett's result \cite[Theorem 1.8]{Ligg75} characterizes the stationary distribution of the half-line open TASEP with step initial condition (denoted by $\eta_{t}^{step}$) as follows:
\begin{theorem}
\label{thm:stationary}
For $\alpha \leq \frac{1}{2}$, the half-line open TASEP with step initial condition $\eta_t^{\mathrm{step}}$ converges to a product measure of i.i.d. Bernoulli random variables with parameter $\alpha$. 
\end{theorem}

When $\alpha>\frac12$, the limiting stationary measure is no longer a product measure. 
Liggett's ergodic theorem \cite[Theorem 1.8]{Ligg75} implies that, starting from the step initial condition, the half-line open TASEP converges weakly to the stationary measure denoted by $\mu_{1/2}^{\alpha}$, which has spatial correlations and is asymptotically product with density $\frac12$ at infinity. 

The measure $\mu_{1/2}^{\alpha}$ admits a matrix-product representation. 
Grosskinsky \cite[Theorem 3.2]{Gro04}, following the DEHP matrix-product ansatz \cite{DEHP93}, constructed stationary measures for the half-line open TASEP with current $c$. 
In the notation of \cite{DMZ18}, these matrix-product measures are denoted by $\bar\nu_c^\alpha$. 
Moreover, \cite[Proposition 3.2]{DMZ18} identifies Grosskinsky's matrix-product measure with Liggett's limiting stationary measure: if $\alpha\geq \frac12$, $\varrho\geq\frac12$, and $c=\varrho(1-\varrho)$, then
\[
\bar\nu_c^\alpha=\mu_\varrho^\alpha,
\]
where $\varrho$ is the stationary particle density.
Taking $\varrho=\frac12$ gives $c=\frac14$. 
Therefore, for the step initial condition and $\alpha>\frac12$, the limiting measure $\mu_{1/2}^{\alpha}$ is precisely the matrix-product measure with current $c=\frac14$.

For completeness, we record the matrix-product form of the limiting measure $\mu_{1/2}^{\alpha}$, which combines Liggett's convergence theorem \cite[Theorem 1.8]{Ligg75} with the matrix-product characterization of Grosskinsky \cite[Theorem 3.2]{Gro04} and the identification in \cite[Proposition 3.2]{DMZ18}.

\begin{theorem}[Liggett; Grosskinsky; Duhart--Mörters--Zimmer]
\label{thm:MPA}
Let $\alpha>\frac12$ and let $c=\frac14$. 
Let $D$, $E$ be non-negative matrices, and let $\bra{w}$, $\ket{v}$ be vectors satisfying
\begin{align}
\label{eq:stationary1}
&DE=c(D+E),\\
\label{eq:stationary2}
&\alpha \bra{w}E=c\bra{w},\\
\label{eq:stationary3}
&(D+E)\ket{v}=\ket{v},\\
\label{eq:stationary4}
&\braket{w}{v}=1.
\end{align}
Then the limiting stationary measure of the half-line open TASEP with step initial condition is the probability measure $\mu_{1/2}^{\alpha}$ determined by the cylinder probabilities
\begin{align}
\label{eq:probab_mu}
\mu_{1/2}^{\alpha}\{\xi \in \{0,1\}^{\mathbb{N}}: \xi_1=\eta_1, \cdots, \xi_L=\eta_L\}
=
\bra{w}\prod_{x=1}^{L}\left(\eta_{x}D+(1-\eta_x)E\right)\ket{v},
\end{align} 
for any $L \in \mathbb{N}$ and $(\eta_1,\cdots,\eta_L)\in\{0,1\}^{L}$.
\end{theorem}

The next result gives an explicit combinatorial evaluation of the cylinder probabilities in the matrix-product formula \eqref{eq:probab_mu}. 
Together with Theorem \ref{thm:MPA}, this provides the large-time cylinder probabilities of the half-line open TASEP with step initial condition in the maximum-current phase $\alpha>\frac12$.
\begin{theorem}
\label{thm:combinatorial}
For any $L\in\mathbb{N}$ and any state $\eta=(\eta_1,\cdots,\eta_L) \in \{0,1\}^{L}$, we define the DEHP-tree of $\eta$ by Definition \ref{def:DEHP-tree}. The stationary measure in \eqref{eq:probab_mu} has the explicit combinatorial formula:
\begin{align}
\label{eq:stationary}
\mu_{1/2}^{\alpha}\{\xi\in\{0,1\}^{\mathbb{N}}: \xi_1=\eta_1, \cdots, \xi_L=\eta_L\}=\sum_{k=0}^{\Psi_{1}(\eta)} \frac{Z_k(\eta)}{2^k} \cdot \left(1+\frac{\alpha-\frac12}{\alpha}k\right),
\end{align}
where $\Psi_{1}(\eta)$ and $Z_k(\eta)$ are defined in Section \ref{ssec:constant-scaling}.
\end{theorem}
In the rest of this section, we will analyze the structure of the matrix product ansatz \eqref{eq:probab_mu} and give a proof of Theorem \ref{thm:combinatorial}.

With rules \eqref{eq:stationary1}, \eqref{eq:stationary2} and \eqref{eq:stationary4}, the probability of any cylinder set $(\eta_1, \cdots, \eta_{L})$ can be directly written as a linear combination of $\bra{w}D^{k}\ket{v}$:
\begin{align*}
\mu_{1/2}^{\alpha}\{\xi\in\{0,1\}^{\mathbb{N}}:\xi_1=\eta_1, \cdots, \xi_L=\eta_L\}=\sum_{k=0}^{\ell(\eta)}d_k(\eta)\bra{w}D^{k}\ket{v},
\end{align*}
where $\ell(\eta)=\#\{i: \eta_i=1\}$. Using \eqref{eq:stationary3}, Grosskinsky \cite{Gro04} has noticed that $\bra{w}D^{k}\ket{v}$ satisfies the following recursion:
\begin{align*}
\bra{w}D^{k+2}\ket{v}=\bra{w}D^{k+1}(D+E)\ket{v}-\bra{w}D^{k}DE\ket{v}=\bra{w}D^{k+1}\ket{v}-c\bra{w}D^{k}\ket{v},
\end{align*}
with initial values $\bra{w}D^{0}\ket{v}=1$ and $\bra{w}D\ket{v}=1-c/\alpha$, and the solutions of this linear recursion are given by the following formulas:
\begin{align}
\label{eq:D_k}
\bra{w}D^{k}\ket{v}=\frac{1}{2^k}\left(1+\frac{\alpha-\frac12}{\alpha}k\right), \quad 0 \leq k \leq \ell(\eta). 
\end{align}
In this section, we are focused on the combinatorial expression of the coefficients $d_k(\eta)$ for any $\eta$ and $0 \leq k \leq \ell(\eta)$.

The DEHP-tree is a graphical way to compute these coefficients $d_k(\eta)$. 
The basic idea is to reduce the non-commutative word in $D$ and $E$ to a normal form containing only powers of $D$ next to the right vector $\ket{v}$. 
This reduction uses two relations. 
First, whenever an $E$ lies to the right of a $D$, the relation
\[
DE=c(D+E)
\]
replaces the adjacent pair $DE$ by two possible terms. 
Second, when an $E$ reaches the left boundary vector, the relation
\[
\alpha\bra{w}E=c\bra{w}
\]
removes this $E$. 
Thus each complete reduction path contributes a product of factors $c$ and $c/\alpha$, and the sum of the weights of all reduction paths ending with $D^k$ gives the coefficient of $\bra{w}D^k\ket{v}$.

We suppose that $\eta=(\eta_1, \cdots, \eta_{L})$ consists of $n$ clusters of $(1,\cdots,1)$ with length $\sigma_i$ and $n+1$ clusters of $(0,\cdots,0)$ with length $\tau_j$, which has the following form:
\begin{align}
\label{eq:etaL}
\eta=(\underbrace{0,\cdots,0}_{\tau_0}, \underbrace{1,\cdots,1}_{\sigma_1}, \underbrace{0,\cdots,0}_{\tau_1}, \underbrace{1,\cdots,1}_{\sigma_2}, \cdots,  \underbrace{0,\cdots,0}_{\tau_{n-1}}, \underbrace{1,\cdots,1}_{\sigma_n}, \underbrace{0,\cdots,0}_{\tau_n})
\end{align}
where $0 \leq n \leq \lceil L/2 \rceil$,  $1 \leq \sigma_i \leq L$ for $i \in \{1, 2, \cdots, n\}$, $1 \leq \tau_j \leq L$ for $j \in \{1, \cdots, n-1\}$ and $0 \leq \tau_0, \tau_n \leq L$. 
Under the matrix-product correspondence, a 1 contributes $D$ and a 0 contributes $E$. 
Therefore the cylinder word associated with \eqref{eq:etaL} is
\[
E^{\tau_0}D^{\sigma_1}E^{\tau_1}\cdots D^{\sigma_n}E^{\tau_n}.
\]
The parameters $\sigma_i$ and $\tau_j$ are thus the lengths of the successive $D$-clusters and $E$-clusters. 
Then the expression of \eqref{eq:probab_mu} is the following form:
\begin{align*}
&\mu_{\alpha}^{1/2}\{\xi \in \{0,1\}^{\mathbb{N}}: \xi_1=\eta_1, \cdots, \xi_L=\eta_L\}\\
=&\bra{w}\underbrace{E\ldots E}_{\tau_0} \underbrace{D \ldots D}_{\sigma_1} \underbrace{E \ldots E}_{\tau_1}  \underbrace{D \ldots D}_{\sigma_2}  \cdots \underbrace{E \ldots E}_{\tau_{n-1}}  \underbrace{D \ldots D}_{\sigma_n}  \underbrace{E \ldots E}_{\tau_n}\ket{v}\\
:=& F(\tau_0,\sigma_1, \tau_1, \cdots, \sigma_n, \tau_n),
\end{align*}
Our main task is to express $F(\tau_0, \sigma_1, \tau_1, \cdots, \sigma_n, \tau_n)$ as a linear combination of $\bra{w}D^{k}\ket{v}$.

We first treat the elementary case of one $D$-block followed by one $E$-block. 
This is the basic local reduction step that will later be iterated over the clusters of a general word. 
For $\sigma_1,\tau_1\geq0$, define
\[
F_1(\sigma_1,\tau_1)
=
\bra{w}D^{\sigma_1}E^{\tau_1}\ket{v}.
\]
The one-layer DEHP-tree records all possible ways of reducing $D^{\sigma_1}E^{\tau_1}$ by repeated use of $DE=c(D+E)$ and $\alpha\bra{w}E=c\bra{w}$.

\begin{figure}
\begin{center}
\begin{tikzpicture}[scale=1]
\draw (0,10) node {$(1,2)$};
\draw[black,line width=1pt,->] (-0.1,9.8) -- (-1,9.3) node[midway, anchor=west]{$c$};
\draw (-1,9) node {$(1,1)$};
\draw[black,line width=1pt,->] (0.1,9.8) -- (1,9.3) node[midway, anchor=east]{$c$};
\draw (1,9) node {$(0,2)$};
\draw[black,line width=1pt,->] (-1.1,8.8) -- (-2,8.3) node[midway, anchor=west]{$c$};
\draw (-2,8) node {$(1,0)$};
\draw[black,line width=1pt,->] (-0.9,8.8) -- (0,8.3) node[midway, anchor=east]{$c$};
\draw (0,8) node {$(0,1)$};
\draw[black,line width=1pt,->] (0.9,8.8) -- (0,8.3) node[midway, anchor=west]{$c/\alpha$};
\draw (-1,7) node {$(0,0)$};
\draw[black,line width=1pt,->] (-0.1,7.8) -- (-1,7.3) node[midway, anchor=west]{$c/\alpha$};
\end{tikzpicture}
\end{center}
\caption{The one-layer DEHP-tree of $DEE$.}
\label{fig:one-layer}
\end{figure}
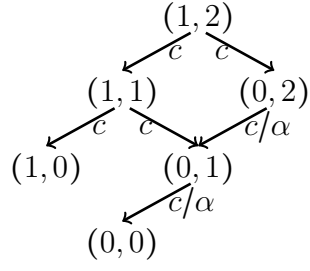

\begin{definition}
\label{def:one-layer}
The {\em  one-layer DEHP-tree} of $(\sigma_1, \tau_1)$ is a rooted tree with {\em  DEHP weight}, which is defined as follows:
\begin{itemize}
\item start from a {\em node} represented by a tuple $(\sigma_1, \tau_1)$, which we call the {\em root};
\item any node $(x,y)$ with $x,y \geq 1$ has two {\em children}: the {\em left child} $(x, y-1)$ and the {\em right child} $(x-1, y)$, the node $(x,y)$ is called parent and we have a {\em weighted directed path} $\swarrow$ (resp. $\searrow$) pointing from the parent to the left child (resp. the right child), where each path has {\em weight} $c$;
\item any node $(0,y)$ with $y \geq 1$ has only one child: the {\em left child} $(0, y-1)$, the parent $(0,y)$ is connected with the left child $(0, y-1)$ by a path weighted by $c/\alpha$;
\item end with the nodes $(x,0)$ where $x=0,1,\cdots,\sigma_1$, which we call the {\em endpoints}.
\end{itemize}
For every node $(x,y)$ reached by this construction, there exists at least one directed path from the root $(\sigma_1,\tau_1)$ to $(x,y)$. We define the weight of such a path as the product of the weights along the path. We define the {\em partition function} of a node $(x,y)$ as the sum of weights of all possible paths from the root $(\sigma_1, \tau_1)$ to that node, and denote it by $Z_{x,y}(\sigma_1, \tau_1)$. In particular, we denote the partition function of endpoints $(k,0)$ as $Z_k(\sigma_1, \tau_1)$.  (See Figure \ref{fig:one-layer} for example). 
\end{definition}

The meaning of the two types of edges is as follows. 
A two-child branching from $(x,y)$ represents the reduction
\[
D^xE^y
=
D^{x-1}(DE)E^{y-1}
=
cD^xE^{y-1}+cD^{x-1}E^y,
\]
and hence gives two children, both with edge weight $c$. 
A one-child branching occurs when no $D$ remains to the left of the unresolved $E$'s. 
Then the boundary relation
\[
\alpha\bra{w}E=c\bra{w}
\]
removes one $E$, giving a single edge with weight $c/\alpha$. 
Consequently, a path ending at $(k,0)$ corresponds to one complete reduction that leaves the monomial $D^k$, and the partition function $Z_k$ is the sum of the weights of all such paths.

\begin{proposition}
\label{prop:one-layer}
We can give the explicit combinatorial formula of $F(\sigma_1, \tau_1)$ by the above one-layer DEHP-tree:
\begin{align}
F(\sigma_1, \tau_1)=\sum_{k=0}^{\sigma_1}Z_k(\sigma_1, \tau_1)\bra{w}D^{k}\ket{v}.
\end{align}
\end{proposition}

\begin{proof}
It is enough to identify the following correspondence between the algebraic reductions and the edges of the one-layer DEHP-tree:
\begin{align*}
(x,y) \iff  \underbrace{D \ldots D}_{x} \underbrace{E \ldots E}_{y} 
\end{align*}
\begin{align*}
\begin{tikzpicture}[scale=1]
\draw (0,1) node {$(x,y)$};
\draw[black,line width=1pt,->] (-0.1,0.8) -- (-1,0.3) node[midway, anchor=west]{$c$};
\draw (-1,0) node {$(x,y-1)$};
\draw[black,line width=1pt,->] (0.1,0.8) -- (1,0.3) node[midway, anchor=east]{$c$};
\draw (1,0) node {$(x-1,y)$};
\draw (2.5,0.5) node {$\iff$};
\draw (3.5,0.5) node {$\eqref{eq:stationary1}$};
\draw (5.5,0.5) node {$\text{for} \ x,y \geq 1$};
\end{tikzpicture}
\end{align*}
\begin{align*}
\begin{tikzpicture}[scale=1]
\draw (0,1) node {$(0,y)$};
\draw[black,line width=1pt,->] (-0.1,0.8) -- (-1,0.3) node[midway, anchor=west]{$c/\alpha$};
\draw (-1,0) node {$(0,y-1)$};
\draw (2.5,0.5) node {$\iff$};
\draw (3.5,0.5) node {$\eqref{eq:stationary2}$};
\draw (5.5,0.5) node {$\text{for} \ y \geq 1$};
\end{tikzpicture}
\end{align*}
\end{proof}

We now iterate the one-layer construction through the successive clusters of the word
\[
E^{\tau_0}D^{\sigma_1}E^{\tau_1}\cdots D^{\sigma_n}E^{\tau_n}.
\]
After processing the first $i$ particle clusters and the corresponding zero clusters, the partial product is reduced to a linear combination of $\bra{w}D^k$. 
The multi-layer DEHP-tree records this recursion cluster by cluster.
We introduce the notation $\tilde{\Psi}_i=\sum_{j=0}^{i}\sigma_j$ and $\tilde{\Phi}_i=\sum_{j=0}^{i}\tau_j$. 
The cumulative quantities $\tilde\Psi_i$ and $\tilde\Phi_i$ record how many $D$'s and $E$'s have been included after the first $i$ particle clusters. 
They are used to locate the layer of the tree corresponding to the first $i$ pairs of $D$- and $E$-blocks.

\begin{figure}
\begin{center}
\begin{tikzpicture}[scale=1]
\draw [black] (-9,10) node {$(3,3)$};
\draw[black,line width=1pt,->] (-9.1,9.8) -- (-10,9.3) node[midway, anchor=west]{$c$};
\draw (-10,9) node {$(3,2)$};
\draw[black,line width=1pt,->] (-8.9,9.8) -- (-8,9.3) node[midway, anchor=east]{$c$};
\draw [black] (-8,9) node {$(2,3)$};
\draw[black,line width=1pt,->] (-10.1,8.8) -- (-11,8.3) node[midway, anchor=west]{$c$};
\draw [black] (-11,8) node {$(3,1)$};
\draw[black,line width=1pt,->] (-9.9,8.8) -- (-9,8.3) node[midway, anchor=east]{$c$};
\draw [black] (-9,8) node {$(2,2)$};
\draw[black,line width=1pt,->] (-11.1,7.8) -- (-12,7.3) node[midway, anchor=west]{$c$};
\draw [black] (-12,7) node {$(3,0)$};
\draw[black,line width=1pt,->] (-10.9,7.8) -- (-10,7.3) node[midway, anchor=east]{$c$};
\draw [black] (-10,7) node {$(2,1)$};
\draw [black] (-11,6) node {$(2,0)$};
\draw[black,line width=1pt,->] (-10.1,6.8) -- (-11,6.3) node[midway, anchor=west]{$c$};
\draw[black,line width=1pt,->] (-9.9,6.8) -- (-9,6.3) node[midway, anchor=east]{$c$};
\draw (-9,6) node {$(1,1)$};
\draw[black,line width=1pt,->] (-9.1,5.8) -- (-10,5.3) node[midway, anchor=west]{$c$};
\draw[black,line width=1pt,->] (-8.9,5.8) -- (-8.1,5.3) node[midway, anchor=east]{$c$};
\draw [black] (-10,5) node {$(1,0)$};
\draw [black] (-8,5) node {$(0,1)$};
\draw [black] (-9,4) node {$(0,0)$};
\draw[black,line width=1pt,->] (-8.1,8.8) -- (-8.9,8.3) node[midway, anchor=west]{$c/\alpha$};
\draw[black,line width=1pt,->] (-9.1,7.8) -- (-9.9,7.3) node[midway, anchor=west]{$c/\alpha$};
\draw[black,line width=1pt,->] (-8.1,4.8) -- (-9,4.3) node[midway, anchor=west]{$c/\alpha$};
\draw [black] (0,10) node {$(5,5)$};
\draw[black,line width=1pt,->] (-0.1,9.8) -- (-1,9.3) node[midway, anchor=west]{$c$};
\draw (-1,9) node {$(5,4)$};
\draw[black,line width=1pt,->] (0.1,9.8) -- (1,9.3) node[midway, anchor=east]{$c$};
\draw [black] (1,9) node {$(4,5)$};
\draw[black,line width=1pt,->] (-1.1,8.8) -- (-2,8.3) node[midway, anchor=west]{$c$};
\draw [black] (-2,8) node {$(5,3)$};
\draw[black,line width=1pt,->] (-0.9,8.8) -- (0,8.3) node[midway, anchor=east]{$c$};
\draw [black] (0,8) node {$(4,4)$};
\draw[black,line width=1pt,->] (-2.1,7.8) -- (-3,7.3) node[midway, anchor=west]{$c$};
\draw [black] (-3,7) node {$(5,2)$};
\draw[black,line width=1pt,->] (-1.9,7.8) -- (-1,7.3) node[midway, anchor=east]{$c$};
\draw [black] (-1,7) node {$(4,3)$};
\draw[black,line width=1pt,->] (-3.1,6.8) -- (-4,6.3) node[midway, anchor=west]{$c$};
\draw (-4,6) node {$(5,1)$};
\draw[black,line width=1pt,->] (-2.9,6.8) -- (-2,6.3) node[midway, anchor=east]{$c$};
\draw [black] (-2,6) node {$(4,2)$};
\draw[black,line width=1pt,->] (-1.1,6.8) -- (-2,6.3) node[midway, anchor=west]{$c$};
\draw[black,line width=1pt,->] (-0.9,6.8) -- (0,6.3) node[midway, anchor=east]{$c$};
\draw (0,6) node {$(3,3)$};
\draw[black,line width=1pt,->] (-4.1,5.8) -- (-5,5.3) node[midway, anchor=west]{$c$};
\draw[black,line width=1pt,->] (-3.9,5.8) -- (-3,5.3) node[midway, anchor=east]{$c$};
\draw[black,line width=1pt,->] (-2.1,5.8) -- (-3,5.3) node[midway, anchor=west]{$c$};
\draw[black,line width=1pt,->] (-1.9,5.8) -- (-1,5.3) node[midway, anchor=east]{$c$};
\draw[black,line width=1pt,->] (-0.1,5.8) -- (-1,5.3) node[midway, anchor=west]{$c$};
\draw[black,line width=1pt,->] (0.1,5.8) -- (0.9,5.3) node[midway, anchor=east]{$c$};
\draw [black] (-5,5) node {$(5,0)$};
\draw (-3,5) node {$(4,1)$};
\draw [black] (-1,5) node {$(3,2)$};
\draw [black] (1,5) node {$(2,3)$};
\draw[black,line width=1pt,->] (-3.1,4.8) -- (-4,4.3) node[midway, anchor=west]{$c$};
\draw[black,line width=1pt,->] (-2.9,4.8) -- (-2,4.3) node[midway, anchor=east]{$c$};
\draw[black,line width=1pt,->] (-1.1,4.8) -- (-2,4.3) node[midway, anchor=west]{$c$};
\draw[black,line width=1pt,->] (-0.9,4.8) -- (0,4.3) node[midway, anchor=east]{$c$};
\draw [black] (-4,4) node {$(4,0)$};
\draw (-2,4) node {$(3,1)$};
\draw [black] (0,4) node {$(2,2)$};
\draw[black,line width=1pt,->] (-2.1,3.8) -- (-3,3.3) node[midway, anchor=west]{$c$};
\draw[black,line width=1pt,->] (-1.9,3.8) -- (-1,3.3) node[midway, anchor=east]{$c$};
\draw[black,line width=1pt,->] (-0.1,3.8) -- (-1,3.3) node[midway, anchor=west]{$c$};
\draw[black,line width=1pt,->] (0.1,3.8) -- (1,3.3) node[midway, anchor=east]{$c$};
\draw [black] (-3,3) node {$(3,0)$};
\draw (-1,3) node {$(2,1)$};
\draw (1,3) node {$(1,2)$};
\draw[black,line width=1pt,->] (-1.1,2.8) -- (-2,2.3) node[midway, anchor=west]{$c$};
\draw[black,line width=1pt,->] (-0.9,2.8) -- (0,2.3) node[midway, anchor=east]{$c$};
\draw[black,line width=1pt,->] (0.9,2.8) -- (0,2.3) node[midway, anchor=west]{$c$};
\draw[black,line width=1pt,->] (1.1,2.8) -- (2,2.3) node[midway, anchor=east]{$c$};
\draw [black] (-2,2) node {$(2,0)$};
\draw (0,2) node {$(1,1)$};
\draw [black] (2,2) node {$(0,2)$};
\draw[black,line width=1pt,->] (-0.1,1.8) -- (-1,1.3) node[midway, anchor=west]{$c$};
\draw[black,line width=1pt,->] (0.1,1.8) -- (1,1.3) node[midway, anchor=east]{$c$};
\draw [black] (-1,1) node {$(1,0)$};
\draw [black] (1,1) node {$(0,1)$};
\draw[black,line width=1pt,->] (0.9,8.8) -- (0.1,8.3) node[midway, anchor=west]{$c/\alpha$};
\draw[black,line width=1pt,->] (-0.1,7.8) -- (-0.9,7.3) node[midway, anchor=west]{$c/\alpha$};
\draw[black,line width=1pt,->] (0.9,4.8) -- (0.1,4.3) node[midway, anchor=west]{$c/\alpha$};
\draw[black,line width=1pt,->] (1.9,1.8) -- (1,1.3) node[midway, anchor=west]{$c/\alpha$};
\draw[black,line width=1pt,->] (0.9,0.8) -- (0,0.3) node[midway, anchor=west]{$c/\alpha$};
\draw [black] (0,0) node {$(0,0)$};
\end{tikzpicture}
\end{center}
\caption{Left panel: The two-layer DEHP-tree of $DEEDDE$; Right panel: The three-layer DEHP-tree of $DEEDDEDDEE$.}
\label{fig:n-layer}
\end{figure}
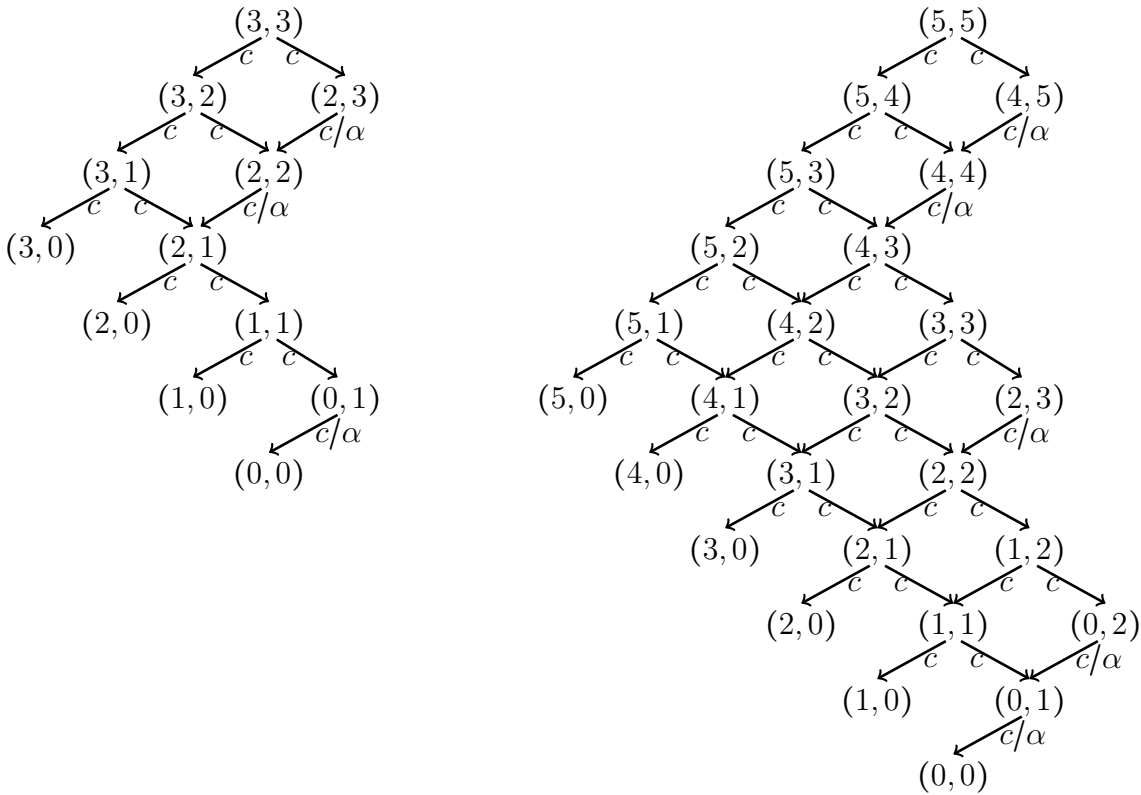

\begin{definition}
\label{def:n-layer tree}
The {\em  multi-layer DEHP-tree} of $(\tilde{\Psi}_n, \tilde{\Phi}_n)$ is a rooted tree with {\em  DEHP weight}, which is defined recursively as follows:
\begin{itemize}
\item start with a {\em one-layer DEHP-tree} of $(\tilde{\Psi}_0, \tilde{\Phi}_0)$;
\item construct the {\em $i$-layer DEHP-tree} of $(\tilde{\Psi}_i, \tilde{\Phi}_i)$ from the {\em $(i-1)$-layer DEHP-tree} of $(\tilde{\Psi}_{i-1}, \tilde{\Phi}_{i-1})$: we add $(\sigma_i, \tau_i)$ to every node of the {\em $(i-1)$-layer DEHP-tree} of $(\tilde{\Psi}_{i-1}, \tilde{\Phi}_{i-1})$ and then continue the tree construction by performing the {\em one-layer DEHP-tree} construction of $(\tilde{\Psi}_i, \tau_i)$.
\item end with the {\em endpoints} $(x,0)$ where $x=0,1,\cdots,\tilde{\Psi}_n$.
\end{itemize}
We define the {\em partition function} of a node $(x,y)$ as the sum of weights of all possible paths from the root $(\tilde{\Psi}_n, \tilde{\Phi}_n)$ to that node, and denote it by $Z_{x,y}(\tilde{\Psi}_n, \tilde{\Phi}_n)$. In particular, we denote the partition function of endpoints $(k,0)$ as $Z_k(\tilde{\Psi}_n, \tilde{\Phi}_n)$.  (See Figure \ref{fig:n-layer} for examples). 
\end{definition}

\begin{proposition}
We can give the explicit formula of $F(\tau_0, \sigma_1, \tau_1, \cdots, \sigma_n, \tau_n)$ by the above {\em  multi-layer DEHP-tree}:
\begin{align}
\label{eq:n-layer}
F(\tau_0, \sigma_1, \tau_1, \cdots, \sigma_n, \tau_n)=\sum_{k=0}^{\tilde{\Psi}_n}Z_k(\tilde{\Psi}_n, \tilde{\Phi}_n)\bra{w}D^{k}\ket{v}.
\end{align}
\end{proposition}

\begin{proof}
By definition, we have $F(\tau_0, \sigma_1, \tau_1, \cdots, \sigma_n, \tau_n)=\bra{w}\prod_{j=0}^{n}\underbrace{D \ldots D}_{\sigma_j} \underbrace{E \ldots E}_{\tau_j}\ket{v}$, where $\sigma_0=0$. We calculate the product step by step. First, $\bra{w} \underbrace{E \ldots E}_{\tau_0}$ is given by the one-layer DEHP-tree: $Z_0(0, \tau_0)\bra{w}$. If we already have $\bra{w}\prod_{j=0}^{i-1}\underbrace{D \ldots D}_{\sigma_j} \underbrace{E \ldots E}_{\tau_j}=\sum_{k=0}^{\tilde{\Psi}_{i-1}}Z_k(\tilde{\Psi}_{i-1}, \tilde{\Phi}_{i-1})\bra{w}D^{k}$, then $\bra{w}\prod_{j=1}^{i}\underbrace{D \ldots D}_{\sigma_j} \underbrace{E \ldots E}_{\tau_j}$ can be expressed by
 \[\sum_{k=0}^{\tilde{\Psi}_{i-1}}Z_k(\tilde{\Psi}_{i-1}, \tilde{\Phi}_{i-1})\bra{w}D^{k}\underbrace{D \ldots D}_{\sigma_i} \underbrace{E \ldots E}_{\tau_i}.\]
By Proposition \ref{prop:one-layer}, applying the one-layer reduction to each summand 
\[
\bra{w}D^{k+\sigma_i}E^{\tau_i}
\]
produces exactly the children and edge weights prescribed in the construction of the $i$-layer DEHP-tree. 
Summing over all incoming coefficients therefore gives
\[
\sum_{k=0}^{\tilde{\Psi}_i}Z_k(\tilde{\Psi}_i,\tilde{\Phi}_i)\bra{w}D^k.
\]
\end{proof}

\begin{corollary}
\label{eq:alpha12}
For any $\eta$ of the form \eqref{eq:etaL}, the stationary measure \eqref{eq:probab_mu} has the explicit formula:
\begin{align*}
\mu_{1/2}^{\alpha}\{\xi\in\{0,1\}^{\mathbb{N}}: \xi_1=\eta_1, \cdots, \xi_L=\eta_L\}=\sum_{k=0}^{\tilde{\Psi}_n} \frac{Z_k(\tilde{\Psi}_n, \tilde{\Phi}_n)}{2^k }\left(1+\frac{\alpha-\frac12}{\alpha}k\right) ,
\end{align*}
where $Z_k(\tilde{\Psi}_n, \tilde{\Phi}_n)$ is given by the partition function of $(k,0)$ in the {\em multi-layer DEHP-tree} of $(\tilde{\Psi}_n, \tilde{\Phi}_n)$.
\end{corollary}

Now we can provide the proof of Theorem \ref{thm:combinatorial}.
\begin{proof}
The preceding multi-layer construction gives the normal-form expansion of the matrix product associated with a cylinder configuration. 
We now identify it with the compact DEHP-tree of Definition \ref{def:DEHP-tree}.

Let $\eta$ be written in the cluster form \eqref{eq:etaL}. 
The root of the multi-layer DEHP-tree is
\[
(\tilde\Psi_n,\tilde\Phi_n)
=
(\Psi_1(\eta),\Phi_0(\eta)),
\]
which is the root in Definition \ref{def:DEHP-tree}. 
During the reduction of the $i$-th zero cluster, the number of remaining $D$'s is equal to $\Psi_i(\eta)$ along the vertical part of the tree, while the number of remaining $E$'s decreases from $\Phi_{i-1}(\eta)$ to $\Phi_i(\eta)$. 
These are precisely the nodes
\[
(x,y)=\bigl(\Psi_i(\eta),y\bigr),
\qquad
\Phi_i(\eta)<y\leq \Phi_{i-1}(\eta),
\]
which have one child in Definition \ref{def:DEHP-tree}; algebraically this is the boundary reduction
\[
\alpha\bra{w}E=c\bra{w}.
\]

All other nodes in the $i$-th layer correspond to unresolved pairs $DE$. 
They therefore have two children, representing the two terms in
\[
DE=c(D+E).
\]
This gives exactly the two-child rule in Definition \ref{def:DEHP-tree}. 
The edge weights also agree: the one-child edges have weight $c/\alpha$, while the two-child edges have weight $c$.

Therefore the compact DEHP-tree of Definition \ref{def:DEHP-tree} is the same weighted tree as the multi-layer DEHP-tree constructed in Definition \ref{def:n-layer tree}, written in the notation of $\eta$. 
By Corollary \ref{eq:alpha12}, the coefficient of $\bra{w}D^k\ket{v}$ is the partition function $Z_k(\eta)$ of paths ending at $(k,0)$. 
Finally, substituting the explicit value
\[
\bra{w}D^k\ket{v}
=
\frac{1}{2^k}\left(1+\frac{\alpha-\frac12}{\alpha}k\right)
\]
from \eqref{eq:D_k} gives \eqref{eq:stationary}.
\end{proof}

\subsection{Explicit formulas for asymptotic distribution of the second-class particles}
\label{ssec:constant-formulas}
Our first main result is the asymptotic behavior of the second-class particles in $\eta_t$ and $\tilde{\eta}_{t}$ when $\alpha\leq \frac12$:
\begin{theorem}
\label{thm:one-shock-stationary-4.1}
Let $\alpha\leq \frac12$. The asymptotic distribution  of the unique second-class particle existing in $\eta_t$ is given by:
\begin{align}
\label{eq:one-shock-exist-1}
\lim_{t\to\infty}\mathbb{P}(f(t)>0)=\sum_{n=0}^{M_1}\binom{M_1+M_2+1}{n}\alpha^n(1-\alpha)^{M_1+M_2+1-n}.
\end{align}

Let $0 \leq m \leq M_2$, the asymptotic distribution  of the height function of second-class particles in $\tilde{\eta}_{t}$ is given by:
\begin{align}
\label{eq:one-shock-height-1}
\lim_{t\to\infty}\mathbb{P}(\mathbf{N}_{2} (1,t)=m)=\binom{M_1+M_2+1}{M_1+M_2+1-m}\alpha^{M_1+M_2+1-m}(1-\alpha)^{m}.
\end{align}
\end{theorem}

\begin{proof}
The proof follows directly from Theorems \ref{thm:tricolor-disappear}, \ref{thm:tricolor-height}, and \ref{thm:stationary}.
\end{proof}

Our second main result is the asymptotic behavior of the second-class particles in $\eta_t$ and $\tilde{\eta}_{t}$ when $\alpha > \frac12$:

\begin{theorem}
\label{thm:one-shock-stationary-4.2}
Let $\alpha > \frac12$ and we use the notation $B_{n}=\{\eta \in \{0,1\}^{M_1+M_2+1}: \ell(\eta)=n\}$. The asymptotic distribution of the unique second-class particle existing in $\eta_t$ is given by:
\begin{align}
\label{eq:one-shock-exist-2}
\lim_{t \to \infty } \mathbb{P}(f(t)>0)=
\sum_{n=0}^{M_1}\sum_{\eta \in B_{n}}\sum_{k=0}^{\Psi_{1}(\eta)} \frac{Z_k(\eta)}{2^k} \cdot \left(1+\frac{\alpha-\frac12}{\alpha}k\right).
\end{align}

Let $0 \leq m \leq M_2$, the asymptotic distribution  of the height function of second-class particles in $\tilde{\eta}_{t}$ is given by:
\begin{align}
\label{eq:one-shock-height-2}
\lim_{t\to\infty}\mathbb{P}(\mathbf{N}_{2} (1,t)=m)
=\sum_{\eta \in B_{M_1+M_2+1-m}}\sum_{k=0}^{\Psi_{1}(\eta)} \frac{Z_k(\eta)}{2^k} \cdot \left(1+\frac{\alpha-\frac12}{\alpha}k\right).
\end{align}
\end{theorem}

\begin{proof}
The proof follows directly from Theorems \ref{thm:tricolor-disappear}, \ref{thm:tricolor-height}, \ref{thm:MPA}, and \ref{thm:combinatorial}.
\end{proof}

Our third main result is the asymptotic behavior of the second-class particle in $\xi_t$ and the third-class particles in $\tilde{\xi}_{t}$ when $ \alpha \leq \frac12$:
\begin{theorem}
\label{thm:two-shocks-stationary-4.1}
Let $ \alpha \leq \frac12$. The asymptotic distribution of the unique second-class particle existing in the half-line open TASEP $\xi_t$ is given by:
\begin{multline}
\label{eq:two-shocks-exist-1}
\lim_{t\to\infty}\mathbb{P}(g(t)>0)=\sum_{m=0}^{N-1}\sum_{n=0}^{M+N-m}\binom{M+N}{m}\binom{M+N+1}{n}\alpha^{m+n}(1-\alpha)^{2M+2N+1-m-n}\\
+\sum_{m=N}^{M+N}\sum_{n=0}^{M}\binom{M+N}{m}\binom{M+N+1}{n}\alpha^{m+n}(1-\alpha)^{2M+2N+1-m-n}.
\end{multline}

Let $0 \leq s \leq N$, the asymptotic distribution of the height function of third-class particles in $\tilde{\xi}_t$ is given by:
\begin{multline}
\label{eq:two-shocks-height-1}
\lim_{t\to\infty}\mathbb{P}(\mathsf{N}_{3} (1,t)=s)=\sum_{n=1}^{s}\binom{M+N}{N-n}\binom{M+N+1}{M+N+1-s+n}\alpha^{M+2N+1-s}(1-\alpha)^{M+s}\\
+\sum_{m=N}^{M+N}\binom{M+N}{m}\binom{M+N+1}{M+N+1-s}\alpha^{M+N+1+m-s}(1-\alpha)^{M+N-m+s}.
\end{multline}
\end{theorem}

\begin{proof}
The proof follows directly from Theorems \ref{thm:two-shock} and \ref{thm:stationary}.
\end{proof}

Our fourth main result is the asymptotic behavior of the second-class particle in $\xi_t$ and the third-class particles in $\tilde{\xi}_{t}$ when $ \alpha > \frac12$:
\begin{theorem}
\label{thm:two-shocks-stationary-4.2}
Let $ \alpha > \frac12$ and we use the notation $B_{m,n}=\{\eta \in \{0,1\}^{2M+2N+1}: \sum_{i=1}^{M+N}\eta_i=m, \sum_{j=M+N+1}^{2M+2N+1}\eta_j=n\}$. The asymptotic distribution of the unique second-class particle existing in the half-line open TASEP $\xi_t$ is given by:
\begin{multline}
\label{eq:two-shocks-exist-2}
\lim_{t\to\infty}\mathbb{P}(g(t)>0)=\sum_{m=0}^{N-1}\sum_{n=0}^{M+N-m}\sum_{\eta \in B_{m,n}} \sum_{k=0}^{\Psi_{1}(\eta)} \frac{Z_k(\eta)}{2^k} \cdot \left(1+\frac{\alpha-\frac12}{\alpha}k\right)  \\
+\sum_{m=N}^{M+N}\sum_{n=0}^{M}\sum_{\eta \in B_{m,n}}\sum_{k=0}^{\Psi_{1}(\eta)} \frac{Z_k(\eta)}{2^k} \cdot \left(1+\frac{\alpha-\frac12}{\alpha}k\right).
\end{multline}

Let $0 \leq s \leq N$, the asymptotic distribution of the height function of third-class particles in $\tilde{\xi}_t$ is given by:
\begin{multline}
\label{eq:two-shocks-height-2}
\lim_{t\to\infty}\mathbb{P}(\mathsf{N}_{3} (1,t)=s)=\sum_{n=1}^{s}\sum_{\eta \in B_{N-n,M+N+1-s+n}} \sum_{k=0}^{\Psi_{1}(\eta)} \frac{Z_k(\eta)}{2^k} \cdot \left(1+\frac{\alpha-\frac12}{\alpha}k\right)  \\
+\sum_{m=N}^{M+N}\sum_{\eta \in B_{m,M+N+1-s}}\sum_{k=0}^{\Psi_{1}(\eta)} \frac{Z_k(\eta)}{2^k} \cdot \left(1+\frac{\alpha-\frac12}{\alpha}k\right).
\end{multline}
\end{theorem}

\begin{proof}
The proof follows directly from Theorems \ref{thm:two-shock}, \ref{thm:MPA}, and \ref{thm:combinatorial}.
\end{proof}

\section{Asymptotic distribution of second-class particles under KPZ-type scaling}
\label{sec:KPZ-scaling}
In this section, we study the distributions of the second-class particles in the shock initial configurations under KPZ-type scaling. First, we need the known asymptotics of the standard half-line TASEP from the step initial condition.

\begin{definition}[\cite{Rai00}, Section 8]
For a $2\times 2$-matrix valued  skew-symmetric kernel,
$$ K(x,y) = \begin{pmatrix}
\kernel_{11}(x,y) & \kernel_{12}(x,y)\\
\kernel_{21}(x, y) & \kernel_{22}(x,y)
\end{pmatrix},\ \ x,y\in \mathbb{X},$$
we define its Fredholm Pfaffian by the series expansion
\begin{equation}
\Pf\big(J+K\big)_{\mathbb{L}^2(\mathbb{X},\mu)} = 1+\sum_{k=1}^{\infty} \frac{1}{k!}
\int_{\mathbb{X}} \dots \int_{\mathbb{X}}  \Pf\Big( K(x_i, x_j)\Big)_{i,j=1}^k \mathrm{d}\mu^{\otimes k}(x_1 \dots x_k),
\label{eq:defFredholmPfaffian}
\end{equation}
provided the series converges, and we recall that for a skew-symmetric  $2k\times 2k$ matrix $A$, its Pfaffian is defined by
\begin{equation}
 \Pf(A) = \frac{1}{2^k k!} \sum_{\sigma\in\mathcal{S}_{2k}} \sign(\sigma) a_{\sigma(1)\sigma(2)}a_{\sigma(3)\sigma(4)} \dots a_{\sigma(2k-1)\sigma(2k)}.
 \label{def:pfaffian}
\end{equation}
The kernel $J$ is defined by 
$$ J(x,y)  = \begin{cases}
\begin{pmatrix}
0 & 1\\-1 & 0
\end{pmatrix} & \text{ if }x=y, \\
0& \text{ if }x\neq y.
\end{cases}$$
\label{def:FredholmPfaffian}
\end{definition}

\begin{definition}
We define the \emph{half-space ${Airy}_2$ process}
$\mathcal{A}_{\varpi}^{HS}$ via its finite dimensional distributions:
\begin{align}
\mathbb{P}(\mathcal{A}_{\varpi}^{HS}(x_1) \leq \eta_1, \cdots, \mathcal{A}_{\varpi}^{HS}(x_k) \leq \eta_k )= \Pf(J - \kernel^{\rm cross})_{\mathbb{L}^2(\mathbb{D}_k(x_1, \cdots, x_k))},
\end{align}
where $0\leq \eta_1<\cdots<\eta_k$,
\[
\mathbb{D}_k(x_1,\cdots,x_k)
=
\big\{(i,x)\in\{1,\ldots,k\}\times\mathbb{R}:x>x_i\big\},
\]
and the right hand side is the Fredholm Pfaffian
(see Definition \ref{def:FredholmPfaffian}) of the kernel
$\kernel^{\rm cross}$, which depends on $\varpi$ and the $\eta_i$.
The kernel $\kernel^{\rm cross}$ was introduced in
\cite[Section 2.5]{BBCS18a} and can be written as
$$
\kernel^{\rm cross}(i,x;j,y)
=
\Ikernel^{\rm cross}(i,x;j,y)
+
\Rkernel^{\rm cross}(i,x;j,y).
$$
We have
\[
\kernel_{21}^{\rm cross}(i,x;j,y)
=
-\kernel_{12}^{\rm cross}(j,y;i,x),
\]
and
\begin{align*}
\Ikernel_{11}^{\rm cross}(i,x;j,y)
={}&
\frac{1}{(2\I\pi)^2}
\int_{\mathcal{C}_{1}^{\pi/3}}\mathrm{d}z
\int_{\mathcal{C}_{1}^{\pi/3}}\mathrm{d}w
\frac{z+\eta_i-w-\eta_j}{z+w+\eta_i+\eta_j}
\frac{z+\varpi+\eta_i}{z+\eta_i}
\frac{w+\varpi+\eta_j}{w+\eta_j}
e^{z^3/3+w^3/3-xz-yw},
\\
\Ikernel_{12}^{\rm cross}(i,x;j,y)
={}&
\frac{1}{(2\I\pi)^2}
\int_{\mathcal{C}_{a_z}^{\pi/3}}\mathrm{d}z
\int_{\mathcal{C}_{a_w}^{\pi/3}}\mathrm{d}w
\frac{z+\eta_i-w+\eta_j}
{2(z+\eta_i)(z+\eta_i+w-\eta_j)}
\frac{z+\varpi+\eta_i}{-w+\varpi+\eta_j}
e^{z^3/3+w^3/3-xz-yw},
\\
\Ikernel_{22}^{\rm cross}(i,x;j,y)
={}&
\frac{1}{(2\I\pi)^2}
\int_{\mathcal{C}_{b_z}^{\pi/3}}\mathrm{d}z
\int_{\mathcal{C}_{b_w}^{\pi/3}}\mathrm{d}w
\frac{z-\eta_i-w+\eta_j}
{4(z-\eta_i+w-\eta_j)}
\frac{e^{z^3/3+w^3/3-xz-yw}}
{(z-\varpi-\eta_i)(w-\varpi-\eta_j)}.
\end{align*}
The contours in $\Ikernel_{12}^{\rm cross}$ are chosen so that
$a_z>-\eta_i$, $a_z+a_w>\eta_j-\eta_i$ and
$a_w<\varpi+\eta_j$.
The contours in $\Ikernel_{22}^{\rm cross}$ are chosen so that:
if $\varpi\leq 0$, then $b_z>\eta_i$ and $b_w>\eta_j$; if
$\varpi>0$, then
$b_z\in(\eta_i,\eta_i+\varpi)$ and
$b_w\in(\eta_j,\eta_j+\varpi)$.

We have
$\Rkernel_{11}^{\rm cross}(i,x;j,y)=0$, and
$\Rkernel_{12}^{\rm cross}(i,x;j,y)=0$ when $i\geq j$.
When $i<j$,
$$
\Rkernel_{12}^{\rm cross}(i,x;j,y)
=
\frac{-\exp\left(
\frac{
-(\eta_i-\eta_j)^4
+6(x+y)(\eta_i-\eta_j)^2
+3(x-y)^2
}{
12(\eta_i-\eta_j)
}
\right)}
{\sqrt{4\pi(\eta_j-\eta_i)}}.
$$
This may also be written as
$$
\Rkernel_{12}^{\rm cross}(i,x;j,y)
=
-\int_{-\infty}^{+\infty}
\mathrm{d}\lambda\,
e^{-\lambda(\eta_i-\eta_j)}
\Ai(x+\lambda)\Ai(y+\lambda).
$$
The kernel $\Rkernel_{22}^{\rm cross}$ is antisymmetric, and when
$x-\eta_i^2>y-\eta_j^2$ we have
\begin{multline*}
\Rkernel_{22}^{\rm cross}(i,x;j,y)
=
\frac{\mathbf{1}_{\{\varpi\leq 0\}}}{4}
\frac{1}{2\I\pi}
\int_{\mathcal{C}_{c_z}^{\pi/3}}\mathrm{d}z\,
\frac{
\exp\big(
(z+\eta_j)^3/3+(\varpi+\eta_i)^3/3
-y(z+\eta_j)-x(\varpi+\eta_i)
\big)}
{\varpi+z}
\\
-\frac{\mathbf{1}_{\{\varpi\leq 0\}}}{4}
\frac{1}{2\I\pi}
\int_{\mathcal{C}_{c_z}^{\pi/3}}\mathrm{d}z\,
\frac{
\exp\big(
(z+\eta_i)^3/3+(\varpi+\eta_j)^3/3
-x(z+\eta_i)-y(\varpi+\eta_j)
\big)}
{\varpi+z}
\\
-\frac{1}{2}
\frac{1}{2\I\pi}
\int_{\mathcal{C}_{d_z}^{\pi/3}}\mathrm{d}z\,
\frac{
z\exp\big(
(z+\eta_i)^3/3+(-z+\eta_j)^3/3
-x(z+\eta_i)-y(-z+\eta_j)
\big)}
{(\varpi+z)(\varpi-z)}
\\
-\frac{\mathbf{1}_{\{\varpi=0\}}}{4}
\exp\big(
\eta_i^3/3+\eta_j^3/3-\eta_i x-\eta_j y
\big).
\end{multline*}
The contours are chosen so that $c_z<-\varpi$ and:
if $\varpi\neq0$, then
$-|\varpi|<d_z<|\varpi|$; if $\varpi=0$, then $d_z>0$.
In the above, we use the notation $\mathcal{C}_{a}^{\varphi}$ for
the union of two semi-infinite rays departing from
$a\in\mathbb{C}$ with angles $\varphi$ and $-\varphi$; we assume
that the contour is oriented from
$a+\infty e^{-\I\varphi}$ to
$a+\infty e^{+\I\varphi}$.
\end{definition}

\begin{theorem}[\cite{BBCS18a,BBCS18b}]
\label{thm:HS-A2}
Let $\{\mathcal{N}(x,t)\}_{x\in\mathbb{Z}_{>0}}$ be the currents of the half-line open TASEP started from empty initial configuration. 
Assume that $u \in \mathbb{R}$. Set $\alpha = \frac{1+2^{4/3}\varpi t^{-1/3}}{2}$, where $\varpi \in \mathbb{R}$. Then we have 
\begin{align}
\lim_{t\to \infty}\frac{\mathcal{N}(2^{\frac{1}{3}}t^{\frac{2}{3}}u, t)-\frac{t}{4}+2^{-\frac{2}{3}}t^{\frac{2}{3}}u-2^{-\frac{4}{3}}t^{\frac{1}{3}}u^2}{-2^{-\frac{4}{3}}t^{\frac{1}{3}}}=\mathcal{A}_{\varpi}^{HS}(u),
\end{align}
where on the right hand side $\mathcal{A}_{\varpi}^{HS}(x)$ stands for the half-space ${Airy}_2$ process, and the convergence is in the sense of finite-dimensional distributions.
\end{theorem}

\begin{theorem}
\label{thm:one-shock-KPZ-5}
Assume that $M_1=M_2=\lfloor at^{\frac23} \rfloor$, $a\in \mathbb{R}_{>0}$. Set $\alpha = \frac{1+2^{4/3}\varpi t^{-1/3}}{2}$, where $\varpi \in \mathbb{R}$. Then the asymptotic distribution  of the unique second-class particle disappearing from the half-line open TASEP $\eta_t$ is given by:
\begin{align}
\label{eq:one-shock-disappear}
\lim_{t \to \infty}\mathbb{P}(f(t)<0)  = \mathbb{P}\left(\mathcal{A}_{\varpi}^{HS}(0)+2^{4/3}a^2 \leq \mathcal{A}_{\varpi}^{HS}(2^{2/3}a)\right).
\end{align}
The asymptotic distribution  of the height function of second-class particles in $\tilde{\eta}_{t}$ is given by:
\begin{align}
\label{one-shock-kpz}
\lim_{t\to\infty}\frac{\mathbf{N}_{2} (1,t)-at^{2/3}}{t^{1/3}} \stackrel{d}{=} \min \left\{2^{-4/3}\mathcal{A}_{\varpi}^{HS}(0)+a^2 - 2^{-4/3}\mathcal{A}_{\varpi}^{HS}(2^{2/3}a),0 \right\}.
\end{align}
\end{theorem}

\begin{proof}
By Theorem \ref{thm:tricolor-disappear}, the asymptotic distribution  of the unique second-class particle disappearing in the half-line open TASEP $\eta_t$ is given by the following formula:
\begin{align}
\label{prob:1-M}
\lim_{t\to\infty}\mathbb{P}(\mathcal{N}(1,t)-\mathcal{N}(2at^{\frac23},t) \geq at^{\frac23}).
\end{align}
By Theorem \ref{thm:HS-A2}, we have 
\begin{align*}
\mathcal{N}(1,t)=\frac{t}{4}-2^{-\frac{4}{3}}t^{\frac{1}{3}}\mathcal{A}_{\varpi}^{HS}(0)+o(t^{\frac{1}{3}}),\\
\mathcal{N}(2at^{\frac23},t)=\frac{t}{4}-a t^{\frac23}+a^2 t^{\frac13}-2^{-\frac{4}{3}}t^{\frac{1}{3}}\mathcal{A}_{\varpi}^{HS}(2^{\frac23}a)+o(t^{\frac{1}{3}}),
\end{align*}
where $o(t^{\frac{1}{3}})$ means that after dividing by $t^{\frac{1}{3}}$ this term will converge to 0 in probability. Therefore, as $t \to \infty$, the probability in \eqref{prob:1-M} converges to
\begin{align*}
\mathbb{P}\left( 2^{-\frac{4}{3}}\mathcal{A}_{\varpi}^{HS}(2^{\frac23}a)-2^{-\frac{4}{3}}\mathcal{A}_{\varpi}^{HS}(0)-a^2 \geq 0\right),
\end{align*}
this gives \eqref{eq:one-shock-disappear}. Similarly, \eqref{one-shock-kpz} follows from Theorem \ref{thm:tricolor-height}:
\begin{align*}
\lim_{t\to\infty} \mathbf{N}_{2} (1,t) \stackrel{d}{=} \lim_{t\to\infty}\min \left\{2at^{2/3}-\mathcal{N}(1,t)+\mathcal{N}(2at^{\frac23},t),at^{2/3} \right\}.
\end{align*}
\end{proof}

\begin{theorem}
\label{thm:two-shocks-KPZ-5}
Assume that $M=N=\lfloor at^{\frac23} \rfloor$, $a\in \mathbb{R}_{>0}$. Set $\alpha = \frac{1+2^{4/3}\varpi t^{-1/3}}{2}$, where $\varpi \in \mathbb{R}$. Then the asymptotic distribution of the unique second-class particle disappearing from the half-line open TASEP $\xi_t$ is given by:
\begin{multline}
\label{eq:two-shocks-disappear}
\lim_{t \to \infty}\mathbb{P}(g(t)<0)  = \mathbb{P}\left(2^{-\frac{4}{3}}\mathcal{A}_{\varpi}^{HS}(2^{\frac53}a)-2^{-\frac{4}{3}}\mathcal{A}_{\varpi}^{HS}(2^{\frac23}a) - 3a^2\right.\\
\left. + \min\left\{2^{-\frac{4}{3}}\mathcal{A}_{\varpi}^{HS}(2^{\frac23}a)-2^{-\frac{4}{3}}\mathcal{A}_{\varpi}^{HS}(0) - a^2,0 \right\} \geq 0\right).
\end{multline}
The asymptotic distribution of the height function of third-class particles in $\tilde{\xi}_t$ is given by:
\begin{multline}
\label{eq:two-shocks-kpz}
\lim_{t\to\infty}\frac{\mathsf{N}_{3} (1,t)-at^{2/3}}{t^{1/3}} \stackrel{d}{=} \min \left\{0, 2^{-4/3}\mathcal{A}_{\varpi}^{HS}(2^{2/3}a)+3a^2 - 2^{-4/3}\mathcal{A}_{\varpi}^{HS}(2^{5/3}a) \right.\\
\left. +\max\{0, 2^{-4/3}\mathcal{A}_{\varpi}^{HS}(0)+a^2 - 2^{-4/3}\mathcal{A}_{\varpi}^{HS}(2^{2/3}a)\}\right\}.
\end{multline}
\end{theorem}

\begin{proof}
Using \eqref{eq:four-color-disappear}, the asymptotic distribution of the unique second-class particle disappearing from the half-line open TASEP $\xi_t$ is given by:
\begin{align}
\lim_{t\to\infty}\mathbb{P}(\mathcal{N}(2at^{\frac23},t)-\mathcal{N}(4at^{\frac23},t)+\min \{\mathcal{N}(1,t)-\mathcal{N}(2at^{\frac23},t), at^{\frac23}\} \geq 2at^{\frac23}).
\end{align}
By Theorem \ref{thm:HS-A2}, we have 
\begin{align*}
\mathcal{N}(1,t)=\frac{t}{4}-2^{-\frac{4}{3}}t^{\frac{1}{3}}\mathcal{A}_{\varpi}^{HS}(0)+o(t^{\frac{1}{3}}),\\
\mathcal{N}(2at^{\frac23},t)=\frac{t}{4}-a t^{\frac23}+a^2 t^{\frac13}-2^{-\frac{4}{3}}t^{\frac{1}{3}}\mathcal{A}_{\varpi}^{HS}(2^{\frac23}a)+o(t^{\frac{1}{3}}),\\
\mathcal{N}(4at^{\frac23},t)=\frac{t}{4}-2a t^{\frac23}+4a^2 t^{\frac13}-2^{-\frac{4}{3}}t^{\frac{1}{3}}\mathcal{A}_{\varpi}^{HS}(2^{\frac53}a)+o(t^{\frac{1}{3}}).
\end{align*}
Then we have 
\begin{multline*}
\lim_{t \to \infty}\mathbb{P}(g(t)<0)  = \mathbb{P}\left(2^{-\frac{4}{3}}\mathcal{A}_{\varpi}^{HS}(2^{\frac53}a)-2^{-\frac{4}{3}}\mathcal{A}_{\varpi}^{HS}(2^{\frac23}a) - 3a^2\right.\\
\left. + \min\left\{2^{-\frac{4}{3}}\mathcal{A}_{\varpi}^{HS}(2^{\frac23}a)-2^{-\frac{4}{3}}\mathcal{A}_{\varpi}^{HS}(0) - a^2,0 \right\} \geq 0\right).
\end{multline*}
Similarly, \eqref{eq:two-shocks-kpz} follows from \eqref{eq:four-color-remain}:
\begin{multline*}
\lim_{t \to \infty}\mathsf{N}_{3} (1,t)  \stackrel{d}{=}\lim_{t \to \infty} \min \bigl\{a t^{\frac23}, 2a t^{\frac23}-\mathcal{N}(2a t^{\frac23},t)+\mathcal{N}(4a t^{\frac23},t)\bigr.\\
\bigl.+\max\{0, a t^{\frac23}-\mathcal{N}(1,t)+\mathcal{N}(2a t^{\frac23},t)\}\bigr\}.
\end{multline*}
\end{proof}

\begin{remark}
The color--position symmetry approach used in this paper is not limited to the one-shock and two-shock initial configurations considered above. 
In principle, one can introduce more color classes and deterministic sorting words to study initial data with several shocks. 
The resulting finite-time identities should again relate the positions of second-class or higher-class particles to particle-counting functions of the standard half-line open TASEP. 
However, for three or more shocks the sorting intervals may overlap in more complicated ways, and the corresponding events involve several nested particle counts. 
A particularly interesting problem is shock merging, where different second-class particles interact on the same scale. 
We leave the systematic analysis of these multi-shock regimes for future work.
\end{remark}

\end{document}